\documentclass{amsart}
\usepackage{graphicx,amsthm,amssymb,a4wide}
\newtheorem{theorem}{Theorem}[section]
\newtheorem{corollary}[theorem]{Corollary}
\newtheorem{lemma}[theorem]{Lemma}
\newtheorem{proposition}[theorem]{Proposition}
\theoremstyle{definition}
\newtheorem{definition}[theorem]{Definition}
\theoremstyle{remark}
\newtheorem{remark}[theorem]{Remark}
\numberwithin{equation}{section}

\def\cal{\mathcal}
\def\bf{\mathbf}

\begin{document}

\title[Algebraic $K$-theory view on $KK$-theory]{Algebraic $K$-theory
view on $KK$-theory}%
\author{Tamaz Kandelaki}%
\address{Author's address: A. Razmadze Mathematical Institute, Georgian
Academy of Sciences, 1, M. Aleksidze St., 380093, Tbilisi, Georgia}%
\email{Kandel@rmi.acnet.ge}

\thanks{The work was partially supported by the FNRS grant
No 7GEPJ065513.01 and the INTAS grant No 00566}%

\dedicatory{To my teacher Hvedri Inassaridze }%
\begin{abstract}
Let a compact group $G$ act on real or complex $C^*$-algebras $A$
and $B$, with $A$ separable and $B$ $\sigma$-unital. We express
the $G$-equivariant Kasparov groups $KK_n(A,B)$ by algebraic
$K$-groups of a certain additive category.
\end{abstract}
\maketitle

\section*{Introduction}

In noncommutative topology and differential geometry one of the major
interest is finding topological invariants of some class of noncommutative
algebras. One of the useful and powerful tool is Kasparov $KK$theory. So,
its comprehensive studying, by methods of other mathematical theories,
may be considered as interesting problem.

In the article \cite{kan2}, it is compendiously given account of
the calculation of Kasparov's $KK$-groups by algebraic $K$-groups.
In this paper we we should like to describe it in detail.

In very first view on the problem calculation of $KK$-theory by algebraic
$K$-theory seems to be improbable. There are some reasons. On the one
hand, algebraic $K$-theory of an algebra is in general quite hard
to calculate. On the other hand, the problem mentioned above seem
seems impossible from the view of Kasparov $KK$-theory.

The key to solve the problem is to find suitable object which is sensible
for both algebraic $K$-theory and $KK$-theory. In this paper we make accent
on the $C^*$-category $\mathrm{Rep}(A,B)$, where $A$ is separable and $B$
is $\sigma $-unital real or complex $C^*$-algebras with action of a
fix compact second countable group. Stimulation factor is that, it was
proved in the paper \cite{kan1}, for the complex algebras case, the
existence of an isomorphism, up to a dimensional shift, between the
topological $K$-theory of the $C^{*}$-category ${\rm Rep}(A;B)$ and
Kasparov's groups $KK$-groups:
\begin{equation}
\label{top}
    KK_n(A,B)=K_{n+1}^t({\rm Rep}(A,B)).
\end{equation}
As it was pointed out in \cite{kan1}, one offers, both in real and
complex cases, to define
algebraic, as well as topological, bivariant $KK$-theories by the
formulas
\begin{equation}\label{alg}
KK_n^a(A,B)=K_{n+1}^a({\rm Rep}(A,B))\;\;\;\mathrm{and}\;\;\;
KK_n^t(A,B)=K_{n+1}^a({\rm Rep}(A,B)),
\end{equation}
where $K_n^a$ and $K_n^t$ are
variants of algebraic and topological $K$-theories respectively,
which will be considered below.

To solve our problem we compare three family of contravariant functors.
These are
\begin{equation}
\label{families}
\{KK^a_{n+1}(-;B)\}_{n\in Z},\;\;\;\{KK^t_{n+1}(-;B)\}_{n\in Z},\;\;\;
\{KK_n(-;B)\}_{n\in Z}.
\end{equation}
One of the main our result is that all these three families are so called
Cuntz-Bott cohomology theories (see definition and properties in Section
\ref{cbcoho}). By Corollary \ref{ntiso} the problem comparison is reduced
in a fixed dimension, where the solation of problem is not difficult. So,
our main result shows that algebraic and topological $K$-theories of
$\mathrm{Rep}(A;B)$ are essential isomorphic to the Kasparov $KK(A;B)$.

Below we shall describe the content of paper.

In the first section there are some definitions, constructions and properties
of $C^*$-categories and $C^*$-categoroidies. These are the objects which
help us to formulate some working principles and interpretations, used in
the next sections. In particular, here we construct category
$\mathrm{Rep}(A;B)$ that is the main subject of our research. Mainly,
the material in this section is given without proof and we refer
\cite{kan3} for details.

 In the next sections we'll need Higson's homotopy invariance
theorem.  Since the complexity of $C^{*}$-algebras is significant to prove
Lemma 3.1.2 and Theorem 3.11 of \cite{hig3}, we couldn't dissemination
the proof of it for real $C^{*}$-algebras.
The purpose of the Section  \ref{homin} is to show that
it can be deduced from Kasparov's homotopy invariant theorem, both for
real and complex algebras. This fact is enough for our purposes.

The next main fact, used in the solution of our problem, is Cuntz-Bott
periodicity theorem . We have shown it in the weaker conditions
then those are in \cite{cunt}. A family of contravariant functors satisfying
this conditions we'll call Cuntz-Bott cohomology. From
Cuntz-Bott periodicity theorem we deduce very simple but useful principle:
Let $\mathbb{E}$ and $\mathbb{E}'$ are Cuntz-Bott (co)homology theories. If
there exists natural isomorphism $\mu _m:E_m\rightarrow E'_m$ for some
$m\in Z$, then there exists natural isomorphism
$$
\mu _n:E_n\rightarrow E'_n
$$
for all $n\in Z$ (Corollary \ref{natiso}).

Excision properties of KK-theory, which was proved in \cite{cusk},
tells us that $KK$theory is Cuntz-Bott cohomology theory relative to the
first argument, and Cuntz-Bott homology relative to the second argument.

In the next two sections we study an interpretation of algebraic and
topological $K$-theories of $C^*$-categories. According to this
interpretation, we define algebraic and topological
$K$-groups for, so called, $C^*$-categoroids. Our definition is a
modification of some arguments from \cite{bass}, \cite{hig1},
\cite{quil}. Let us consider it more in detail.

Let $a$ and $a'$ be objects in an additive
$C^{*}$-category $A$ and $I$ be a closed ideal. We say that
$a\leq a'$ if there exists a morphism $s:a\rightarrow a'$
such that $s^{*}s=1_a$ (such type morphism is said to be
isometry). Denote by $\mathcal{L}(a)$ (resp. $I(a)$) the
$C^{*}$-algebra $\hom _A(a,a)$ (resp $\hom _I(a,a)$). We have a
correctly defined inductive system of abelian groups
$\{K_{n}^a(L(a)),\sigma _{aa'}\}_a$ and
$\{K_{n}^t(I(a)),\sigma _{aa^{\prime }}\}_a$, where $K_{n}^a$
(resp. $K_{n}^t$) are usual algebraic (resp. topological)
$K$-theory groups of the algebra (resp. $C^{*}$-algebra) ${\cal
L}(a)$. We suppose that
$$
K_{n}^a(A)=\underrightarrow{\lim}_aK_{n}^a({\cal L}(a))\;\;\;
(\mathrm{resp.})\;\;\;
K_{n}^t(A)=\underrightarrow{\lim}_aK_{n}^t({\cal L}(a))).
$$
So defined algebraic $K$-groups are naturally isomorphic to the
Quillen's $K$-groups $K_{n}^Q(A)$ (with respect to the class of
all split short exact sequences), when $n\geq 0$; and the case topological
$K$-theory gives us an interpretation Karoubi's topological $K$-theory,
\cite{kar}. One can generalize this definition for an ideal $I$ which does
not depend on the choice of the enveloping additive $C^{*}$-category of $I$.
These $K$-theories have excision property generalizing the
analogous property of algebraic and topological $K$-theories of
$C^{*}$-algebras, \cite{kar}, \cite{suw}. There is a natural
transformation $\theta _{n}:K_{n}^a\rightarrow K_{n}^t$ which is a
generalization of the classical natural transformation between
algebraic and topological $K$-theories.

The section \ref{auxiso} of the paper is relatively difficult. In this part
we'll show that
$$
\{KK^a_n(-;B)\}_{n\in Z},\;\;\mathrm{and}\;\;\;\{KK^t_{n+1}(-;B)\}_{n\in Z}
$$
have weak excision property.

In the section 7 we'll prove our main result which says that
all three theories \ref{families} are isomorphic Cuntz-Bott cohomology
theories.

\section{$C^*$-categories and $C^*$categoroids}

In this section we give some elementary properties of $C^{*}$-categories
and $C^{*}$ -categoroids, a natural categorical generalization of
unital $C^{*}$ -algebras and $C^{*}$-algebras. We'll mainly give basic
definitions, constructions and properties without proofs here but we indicate
for details the article \cite{kan3}.

Recall that {\it the diagram scheme} $D$ consists a class of objects
$\mathrm{Ob}D$ and a set $\hom (a,b)$ for any $a,b\in \mathrm{Ob}D$.
By a $k$-{\it scheme } we mean a diagram scheme $D$ such that
$\hom (a,b)$ has the structure of a $k$-linear space, where $k=\mathbb{R}$
or $\mathbb{C}$, the fields of real and complex numbers respectively.
$D$ is called an {\it involutive} $k$-scheme if:
\begin{itemize}
\item  an anti-linear map $*:\hom (a,b)\rightarrow \hom (b,a)$ is given for
       each $a,b\in \mathrm{Ob}D$.
\item  the bilinear composition law
\begin{align*}
    \hom (a,b)\times \hom (b,a) & \rightarrow \hom (a,a) \\
    \hom (a,b)\times \hom (b,b) & \rightarrow \hom (a,b) \\
    \hom (a,a)\times \hom (a,b) & \rightarrow \hom (a,b)
\end{align*}
        is associative for any $a,b\in \mathrm{Ob}D$.

\item   $(f^{*})^{*}=f$, and $(fg)^{*}=g^{*}f^{*}$ if the composition $fg$
        exists.
\end{itemize}
\begin{definition}\label{ae}
By a $C^{*}$-{\it scheme} is meant an involutive $k$-scheme $D$ such that:

$1)$ $\hom (a,b)$ is a Banach space;
$2)$ involution is isometry;
$3)$ $\|f\|^2\leq \|f^{*}f+g^{*}g\|^2$ for any $f\in \hom (a,b)$ and
     $g\in \hom (a,b')$ \cite{pal} \cite{kas1};
$4)$ the morphism $f^{*}f$ is a positive element in the $C^{*}$-algebra
     $\hom (a,a)$ for any $f\in \hom (a,b)$ and $a,b\in ObD$.
\end{definition}

A diagram scheme ${\cal D}$ is called {\it a categoroid} if it satisfies
all the axioms of category except the existence of identities of objects.
Let $a$ and $b$ be objects from ${\cal D}$. Then $\hom (a,b)$ denotes
a set of morphisms from $a$ to $b.$ The definition of morphisms between
categoroids is analogous to that of a functor, and is called
{\it a functoroid}. If ${\cal F:D\rightarrow D^{\prime }\ }$ is a functoroid
of categoroids and there exists composition of morphisms $f, g$
in ${\cal D}$, then ${\cal F}(fg)={\cal F}(f){\cal F}(g).$
\begin{definition}\label{af}
A categoroid $A$ is called a $C^{*}$-categoroid if it has the structure of a
$C^{*}$-scheme such that \\ 1) the composition of morphisms
is bilinear and $\Vert fg\Vert \leq \Vert f\Vert \cdot \Vert g\Vert $
if there exists composition of $f$ and $g$; \\ 2) If $A$ is both a category
and a $C^{*}$-categoroid, then it is called a $C^{*}$-category.
\end{definition}

{\it Remark./}
With the term ''$C^{*}$-categoroid'' will
mean a $C^{*}$-categoroid with a small underlying categoroid,
while the term ''large $C^{*}$-categoroid'' will mean a
$C^{*}$-categoroid with a usual underlying categoroid.

$1)$ The category with all Hilbert spaces as objects and all bounded linear maps
as morphisms is a large $C^{*}$-category denoted by ${\cal H}$.

$2)$ Let $A$ be a $C^{*}$-algebra. The category $\mathcal{H}(A)$
with all right $A$-modules as objects and all bounded $A$-linear maps with
adjoint as morphisms is a large $C^{*}$-category.

$3)$ $C^{*}$-algebra is a $C^{*}$-categoroid with one object $\diamond$
and elements of the $C^{*}$-algebra as morphisms.

$4)$ The category of Hilbert spaces as objects and compact linear maps
as morphisms has the structure of a large $C^{*}$-categoroid.

Let $A$ and $B$ be $C^{*}$-categoroids. A $*$-functoroid
$\mathcal{F}:A\rightarrow B$ is given if ${\cal F}$ maps the objects and
morphisms of $A$ into the objects and morphisms of $B$, respectively,
so that:

$\mathrm{a}$) ${\cal F}(fg)={\cal F}(f){\cal F}(g)$;

$\mathrm{b}$) ${\cal F}(f+g)={\cal F}(f)+{\cal F}(g)$;

$\mathrm{c}$) ${\cal F}(\lambda f)=\lambda {\cal F}(f);$

$\mathrm{d}$) ${\cal F}(f^{*})={\cal F}(f)^{*}$ when
the left side is defined.

If $A$ and $B$ are categories and
${\cal F}(1_a)=1_{{\cal F}(1_a)}$ for any $a\in \mathrm{Ob}A$, then
$\mathcal{F}$ is called a $*$-{\rm functor}. We say that a $*$-functoroid
between $C^{*}$-categoroids is faithful if canonical maps between objects and
between morphisms are injections.

Any $*$-functoroid is norm decreasing. Moreover, the faithful $*$-functoroid
is norm preserving.

\subsection{An Ideal and the Smallest Categorization}

Let $A$ be a $C^{*}$-categoroid and $I\subset \mathrm{Mor}A$. Let
$(a,b)_I=\hom (a,b)\cap I$. Then $I$ is called a left ideal if
$(a,b)_I$ is a linear subspace of $\hom (a,b)$ and $f\in
(a,b)_I,g\in \hom (b,c)$ implies $gf\in (a,c)_I$. A right ideal is
defined similarly. $I$ is a two-sided ideal if it is both a left
and a right ideal. An ideal $I$ is closed if $(a,b)_I$ is closed
in $\hom (a,b)$ for each pair of objects. $I$ determines an
equivalence relation on the morphisms of $A:f\sim g$, if $f-g\in
I$. If $I=I^{*}$ is an ideal of $A$, the set of equivalence
classes $A/I$ can be made into a ${*}$-categoroid in a unique way
by requiring that the canonical map $f\mapsto \bar f$ of
$A\rightarrow A/I$ be a ${*}$-functoroid. $A/I$ can be made into a
normed ${*}$-categoroid by defining
$$
|\bar f| =\sup _{g\in \bar f}|g|.
$$

Arguing as for $C^{*}$-algebras, one can show that
if $A$ be a $C^{*}$-categoroid and $I$ a closed two-sided ideal of $A$.
Then $I=I^{*}$ and $A/I$ is a $C^{*}$-categoroid.

An ideal $I$ in $A$ is called {\it essential} if the intersection
$I\cap J\neq 0$ for each nonzero ideal $J\subset A$.

Let $A$ be a $C^{*}$-categoroid. $C^{*}$-category $B$ is called {\it
categorization} of $A$ if $A$ is contained in $B$ as an essential ideal.

There exists {\it the smallest categorization} $A^{+}$ of $A$.
This $C^{*}$-category has the same objects as $A$, while $\hom
_{A^{+}}(a,a)=\hom_A (a,a)^+$ if the $C^{*}$-algebra $\hom (a,a)$
is non-unital and remains otherwise.

\subsection{Realization of a $C^{*}$-Categoroid}

Any $*$-functoroid $\mathcal{F}:A\rightarrow \mathcal{H}(B)$ is called
{\it a representation}, where $\mathcal{H}(B)$ is a large $C^{*}$-category
of right Hilbert $B$-modules over the $C^{*}$-algebra $B$ and
$B$-linear maps which have adjoint. If $\mathcal{F}$ is
faithful, then it is called {\it a faithful representation} or an
$*$-{\it imbedding. }

Let $A$ be a $C^{*}$-categoroid. There exist a $C^{*}$-algebra $\mathcal{A}$
and a faithful representation

\begin{equation}\label{aa}
    {\cal F}:A\rightarrow {\cal H}({\cal A}).
\end{equation}

\subsection{Stable $C^{*}$-Algebra of a $C^{*}$-Categoroid and Its
Representation}

Let $A$ be a $C^{*}$-categoroid and $S^0(A)$ be a $*$-algebra of finite
$I\times I$-matrices, i.e., of matrices ($a_{ij})_{i,j\in I}$ with only
finite nonzero entries, where $a_{ij}:i\rightarrow j$ is a morphism and
$I $ is the set of objects of $A$. As it was pointed above there exists a
faithful $*$-functoroid from $A$ into the $C^{*} $-category of Hilbert
$\mathcal{A}$-modules and bounded $\mathcal{A}$-homomorphisms, i.e., there is
an injection $*$-functoroid on objects and morphisms. Hence there exists
a $*$-monomorphism of $S^0(A)$ into $L_{\mathcal{A}}(\oplus _iE_i)$, where
$L_{\mathcal{A}}(\oplus E_i)$ is an $C^{*}$-algebra of all bounded
$\mathcal{A}$-homomorphisms on $\oplus _iE_i$ which have adjoint. Therefore
$S^0(A)$ has a (unique) $C^{*}$-norm induced by
$L_{\mathcal{A}}(\oplus E_i)$. A completion of $S^0(A)$ gives a
$C^{*}$-algebra $S(A)$ called a {\it stable} $C^{*}$-algebra of
the $C^{*}$-categoroid $A$.

Let $L=\{{a_1,\dots,a_n\}}$ be a finite subset of objects from $A$, and $A_L$
be a full $C^{*}$-sub-categoroid of $A$ with $L$ as a set of objects. If $R$ is
another finite subset of $\mathrm{Ob}A$ and $L\subset R$, then we have
the canonical ${*}$-homomorphism $i_{LR}:S^0(A_L)\rightarrow S^0(A_R)$
defined by $(f_{ij})\mapsto (g_{kl})$, where $g_{kl}=f_{ij}$ for
$(k,l)=(i,j)$, and is zero otherwise. Thus we obtain the direct system
$\{S^0(A_L),i_{LR}\}$of pre-$C^{*}$-algebras. One has the following
statements:

$\mathrm{a}$) Each element from $S(A)$ can be represented as the
$I\times I$-matrix $(a_{kl})$, where $a_{kl}\in \hom (k,l).$

$\mathrm{b}$) Let $L=\{{a_1,\dots,a_n\}%
}$ be a finite subset of objects from $A$, and $A_L$ be a full
$C^{*}$-sub-categoroid of $A$ with $L$ as a set of objects.
Then $S^0(A_L)$ is a $C^{*}$-algebra.

$\mathrm{c}$) A stable $C^{*}$-algebra $S(A)$ is the direct limit of
the direct system $\{S^0(A_L),i_{LR}\}$ of $C^{*}$-algebras.

\subsection{ Multiplier $C^{*}$-Category of a $C^{*}$-Categoroid}

To construct a multiplier category, we need the following definition.

Let $A$ be a $C^{*}$-algebra and $P=\{p_i\}_{i\in I}$ be a set of
projections of $M(A)$. We say that $P$ is {\it a strictly complete set of
projections} if:
\begin{enumerate}
\item[(a)] $P$ is orthogonal, i.e., $p_ip_j=0$ for every $i\neq j$ from $I$;

\item[(b)] the net $\{p_l\}$ converges strictly to $1$, where $l$ is a
    finite subset of $I$ and $p_l=\sum\limits_{a\in l}p_a$;
    This means that $\sum\limits_{a\in I}p_a=1$ in the strict topology of
    $M(A)$.
\end{enumerate}

Let $A$ be a $C^{*}$-categoroid. A $C^{*}$-category $M(A)$ is
called {\it a multiplier }$C^{*}$-{\it category} of $A$ if $A$ is a closed
two-sided ideal in $M(A)$ and has following universal property: let $D$ be a
$C^{*}$-categoroid containing $A$ as a closed two-sided ideal; then there
exists a unique $*$-functoroid $d:D\rightarrow M(A)$ which is the identity
map on $A$, such that the diagram
\begin{equation}\label{mg}
    \begin{array}{ccc}
        A & \subset & M(A) \\
        \cap & \nearrow &  \\
        D &  &
    \end{array}
\end{equation}
is commutative.

Now, we want to construct the multiplier $C^{*}$-category for any
$C^{*}$-categoroid.

Let $A$ be a $C^{*}$-categoroid. There is an orthogonal strictly complete
set of projections $\mathcal{P}=\{p_a\}_{a\in \mathrm{Ob}A}$ in $M(S(A))$.
Indeed, let $A^{+}$ be the smallest categorization of $A.$ Since
$S(A)$ is an essential ideal in $S(A^{+})$ the canonical $*$-homomorphism
$S(A^{+})\rightarrow M(S(A))$ is injective. For each object $a\in A^{+}$
there is the identity morphism $1_a$ which gives a set of orthogonal
projections $\{p_a\}_{a\in \mathrm{Ob}A}\subset S(A^{+})$. (The image of
$1_a$ in $M(S(A))$ is denoted by $p_a)$. The sums
$p_l=\sum\limits_{a\in l}p_a$ are projections for any finite set $l$
of objects in $A$. Thus we have a set
of projections $\{p_l\}$ in $S(A^{+})$ and $M(S(A))$. It can be
easily shown that $\{p_l\}$ is approximate unit for $S(A^{+})$ so that
$\{1_l\}\rightarrow 1$ strictly in $M(S(A^{+}))$. On the other hand, we
have the strictly continuous unital homomorphism
$M(S(A^{+}))\rightarrow M(S(A))$ since $S(A)$ is an essential ideal
in $S(A^{+})$.
Thus $\{p_l\}$ converges strictly to $1$ in $M(S(A))$.

The set of projections $\{p_a\}_{a\in \mathrm{ob}A}$ is called
{\it the standard complete set of projections} in $M(S(A))$. Let $M(A)$
be a $C^{*}$-category with $\mathrm{Ob}A$ as the set of objects and
the set of
elements of the form $p_{a^{\prime }}\cdot f\cdot p_a$ from $M(S(A))$ as
the set of morphisms from $a$ to $a^{\prime }$. The $C^{*}$-structure is
induced by the corresponding structure of $M(S(A))$.

The sub-category consisting elements from $S(A)$ with property
$x=p_{a'}\cdot x\cdot p_a$ form a closed two-sided ideal $\mathcal{A}$ in
$M(A)$ which is naturally isomorphic to $A$.

Let $B$ be a $C^{*}$-category and $A$ be a closed ideal in $B$.
The category $B$ is the multiplier $C^{*}$-category of $A$ if and only if for
any set with two objects $R=\{a,a'\}$ of $\mathrm{Ob}B$ the canonical
$*$-homomorphism $S(B|_R)\approx M(S(A|_R)$ is an isomorphism.

\subsection{Existence of Direct Limits in the Category Categoroids
\label{dirl}}

Let $\mathbb{C}^*$ be the category of $C^*$-categoroids and
$*$-functoroids. Bellow we'll prove existence of direct
limits of $C^*$-categoroids over directed sets.

Let $S$ be a directed set and $\{A_\alpha,\tau_{\alpha\beta}\}$ be an
inductive system of $C^*$-categoroids. Let $\tau_{\alpha\beta}$ be
monomorphisms on morphisms, i.e. $\tau_{\alpha\beta}(a)=0$ iff
$a=0$, for all $\alpha\le\beta$. Let $A^\circ$ be the
categoroid with $\mathrm{Ob}A^\circ=\varinjlim_\alpha\mathrm{Ob}A_\alpha$
and
$\hom_{A^\circ}([A_\alpha],[A'_\alpha])=
\varinjlim_\alpha(A_\alpha,A'_\alpha)$.
It is clear that $A^\circ$ is a $*$-categoroid. Since $\tau_{\alpha\beta}$
are monomorphisms one can define
$$
||[f_\alpha]||=||f_\alpha||
$$
for all $\alpha$, which gives a pre-$C^*$-norm structure on
$\hom_{A^\circ}([A_\alpha],[A'_\alpha])$. So after completion we have
$\hom([A_\alpha],[A'_\alpha])$. One gets in
such a way a $C^*$-categoroid $A$ is obtained. For each
$\alpha$ we have the $*$-functoroid $\tau_\alpha:A_\alpha\rightarrow A$.

Let us now consider the general case. For each $A_\alpha$ let
$I_{\alpha\beta}$ be the $C^*$-ideal $\ker\tau_{\beta\alpha}$ in
$A_\alpha$. It is clear that $I_{\alpha\beta}\subset
I_{\alpha\gamma}$ if $\beta\le\gamma$. Thus $\{I_{\alpha\beta}\}$
is a direct system of ideals and $I_\alpha=\overline{\bigcup_\beta
I_{\alpha\beta}}$ is a two-sided ideal in $A_\alpha$. Put
$A'_\alpha={A_\alpha}/{I_\alpha}$. Since
$I_\alpha={\tau_{\beta\alpha}}^{-1}(I_\beta)$, the induced
$*$-functor $\tau'_{\alpha\beta}:A'_\alpha\rightarrow A'_\beta$ is
a $*$-mono-functor. As above, $\{A'_\alpha,\tau'_{\alpha\beta}\}$
produces $\{A',\tau'_\alpha\}$, which is the direct limit. Let us show that
the pair $\{A',\tau'_\alpha k_\alpha\}$, where
$k_\alpha:A_\alpha\rightarrow A'_\alpha$ is the canonical $*$-functoroid,
is the direct limit of $\{A_\alpha,\tau_{\alpha\beta}\}$. We
have the commutative diagram
$$
\begin{array}{ccc}
  A_\alpha & \stackrel{\tau_{\alpha\beta}}{\rightarrow}&  A_\beta\\
  \downarrow ^{k_\alpha} & & \downarrow ^{k_\beta} \\
  A'_\alpha & \stackrel{\tau'_{\alpha\beta}}{\rightarrow} & A'_\beta
\end{array}
$$
If $\{B,\nu_\alpha\}$ is a direct system such that

(1) $\nu_\alpha: A_\alpha\rightarrow B$ is $*$-functoroid;

(2) $\nu_\alpha=\nu_\beta\tau_{\alpha\beta}$, then
$\ker\nu_\alpha\supset{\tau_{\beta\alpha}}^{-1}(\ker\nu_\beta)\supset
\ker\tau_{\beta\alpha}$, so $\ker\nu_\alpha\supset I_\alpha$.

Hence there are $*$-functoroids $\nu'_\alpha$ such that
$\nu_\alpha=k_\alpha\nu'_\alpha$; on the other hand
$\nu'_\alpha=\tau'_{\beta\alpha}\nu'_\beta$. Thus there is a
canonical $*$-functoroid $\nu:A'\rightarrow B$ such that
$\nu\tau_\alpha=\nu_\alpha$.

 An inductive system of $C^*$-categories $\{A_\alpha,\tau_{\alpha\beta}\}$
 and let $\{A;\tau_\alpha \}_\alpha$ be a set, where
 $\tau_\alpha :A_\alpha\rightarrow A$
 is functor for each $\alpha$. The set $\{A;\tau_\alpha \}_\alpha$
 is said to be

\begin{enumerate}
  \item  {\it a direct quasi-limit} if for each
$\alpha\le\beta$, $\tau_\beta\cdot\tau_{\alpha\beta}=\tau_\alpha$;
  \item {\it filled}, if for any morphism $x\in A$, the morphism
$\tau(x)$ is invertible iff there exists
$\beta\ge\alpha$ such that $\tau_{\alpha\beta}(x)$ is invertible
in $A_\beta$;
  \item {\it full},
iff for every $x\in A$ and $\epsilon>0$ there exist $\alpha$ and
$x_\alpha\in A_\alpha$ such that
$\|\tau_\alpha(x_\alpha)-x\|<\epsilon$.
\end{enumerate}
Let $\{A,\tau_\alpha\}$ be the direct limit of
the inductive system of $C^*$-categories $\{A_\alpha,\tau_{\alpha\beta}\}$.
Then $\{A,\tau_\alpha\}$ is a filled and full direct quasi-limit. Indeed,
fullness comes from the construction of the direct
limit. It is also filled. Indeed, let $a_\alpha$ be morphisms from
$A_\alpha$ such that $\tau_\alpha(a_\alpha)$ are invertible in $A$. Then by
construction of the direct limit there exists $\beta$ such that
$k_\beta(\tau_{\alpha\beta}(a_\alpha))$ is invertible in
$A'_\beta$, hence there exists an element $a_\beta$ such that
$k_\beta(a_\beta)$ is inverse of
$k_\beta(\tau_{\alpha\beta}(a_\alpha))$. Thus
$\tau_{\alpha\beta}(a_\alpha)\cdot a_\beta=1+i_\beta$,
$a_\beta\cdot\tau_{\alpha\beta}(a_\alpha)=1+i'_\beta$, where
$i_\beta,i'_\beta\in I_\beta$. In general $\|i_\beta\|<1$,
$\|i'_\beta\|<1$ for enough large $\beta$ since
$I_\beta=\overline{\bigcup_\gamma I_{\beta\gamma}}$.
Therefore $\tau_{\alpha\beta}(a_\alpha)\cdot
a_\beta$ and $a_\beta\cdot\tau_{\alpha\beta}(a_\alpha)$ are
invertible elements, so $\tau_{\alpha\beta}(a_\alpha)$ is
invertible.

\subsection{Additive and Pseudo-abelian $C^*$-categories and
$C^*$-categoroids}

A $C^{*}$-categoroid $A$ is said {\it additive }$C^{*}$-{\it categoroid},
if there exists an additive $C^*$-categoroid containing $A$ as a closed
two-sided ideal. Of course, in this situation the multiplier $C^*$-category
must be additive $C^*$-category. A functoroid $f:A\rightarrow A'$ is said
additive if it is restriction of some additive functor between additive $C^*$-
categories containing respectively $A$ and $A'$ as ideals.

Let $a'$ and $a$ be objects in $A$ then each element in
$\hom(a\oplus a')$ may be uniquely represented by matrix of form
$$
\begin{pmatrix}
  l_{aa} & l_{a'a} \\
  l_{aa'} & l_{a'a'}
\end{pmatrix}
$$
where $l_{aa}\in \hom(a,a)$, $l_{aa'}\in \hom(a,a')$, $l_{a'a}\in
\hom(a',a)$, $l_{aa}\in \hom(a',a')$.

\begin{lemma}
\label{print} Let $A$ and $B$ be additive $C^*$-categoroids and
$f:A\rightarrow B$ be an additive $*$-functoroid. Then $f$ is
$*$-isomorphism if and only if $f$ is bijection on the objects and
induced $*$-homomorphism of $C^*$-algebras
$f_a:{\cal L}(a)\rightarrow {\cal L}(f(a))$ is an $*$-isomorphism
for any object $a$ in $A$.
\end{lemma}

\begin{proof}
We only must show that this conditions is sufficient.
It is enough to show that for any pair objects $a,a'$ in $A$, the
linear map
$$
f_{aa'}:\hom(a,a')\rightarrow \hom (f(a),f(a'))
$$
is an isomorphism. i. if $f_{aa'}(\alpha )=0$ then $f_{aa'}(\alpha
)f_{aa'}(\alpha ^*)=f_{aa'}(\alpha ^*)f_{aa'}(\alpha )=0$. This
implies that $\alpha ^*\alpha=\alpha \alpha ^*=0$. Therefore
$\alpha =0$. This means that $f_{aa'}$ is injective. ii. Let
$\beta \in \hom (a,a')$. Consider the matrix
$$
\beta '=
\begin{pmatrix}
  0 & 0 \\
  \beta & 0
\end{pmatrix}
$$
the element of $\mathcal{L}(f(a)\oplus f(a')$. Since
$f_a:\mathcal{L}(a\oplus a')\rightarrow \mathcal{L}(f(a)\oplus f(a'))$
is $*$-isomorphism there exists unique $\alpha \in $ such that
$f_{a}(\alpha)=\beta '$ where $\alpha \in \mathcal{L}(a\oplus a')$.
Let
$$
\alpha =
\begin{pmatrix}
 \alpha _{aa} & \alpha _{a'a} \\
  \alpha _{aa'} & \alpha _{a'a'}
\end{pmatrix}
$$
Since $f_{a\oplus a'}$ is $*$-isomorphism and $\beta $ has unique
representation in matrix form, $\alpha $ must has the form
$$
\alpha =
\begin{pmatrix}
 0 & 0 \\
  \alpha _{aa'} & 0
\end{pmatrix}
$$
Therefore $f_{aa'}(\alpha _{aa'})=\beta $
\end{proof}

We say that a $C^{*}$-category $A$ {\it has enough isometries} if
for any projection $p\in \mathcal{L}(a)$ and any object $a$ in $A$
there exists isometry $s:a'\rightarrow a$ satisfying
equality $ss^{*}=p$. If an additive $C^{*}$-category has enough
isometries then it will be said {\it pseudo-abelian} $C^{*}$-category.

The universal pseudo-additive $C^{*}$-category $\mathrm{P}(A)$ of an
additive $C^{*}$-category $A$ may be constructed by the following
manner \cite{kan1} :

i. objects are all pairs of form $(a,p)$, $a\in \mathrm{Ob}A$,
where $p\in \mathcal{L}(a)$, such that $p^2=p$ and $p^{*}=p$;

ii. We say that triple $(f,p',p)$ (or simply, $f$) is a morphism
in ${\rm P}(A)$ from $(a',p')$ into $(a,p)$ if $f:a'\rightarrow a$
is an morphism in $A$ and $fp'=pf=f$.\\The direct sum of
objects and morphisms are defined by the formulas:
$$
(a,p)\oplus (a',p')=(a\oplus a',p\oplus p')\mathrm{and}
(f,p,q)\oplus (f',p',q')=(f\oplus f',p\oplus p',q\oplus q').
$$
The category $\mathrm{P}(A)$ has natural structure of an additive
C$^{*}$-category inherited from the $C^{*}$-category $A$; Besides,
$A$ can be identified with a full subcategory of $\mathrm{P}(A)$ by
the map of objects $a\mapsto (a,1)$. Since
$(q,q,p):(a,q)\rightarrow (a,p)$ is isometry with
$(q,q,p)(q,q,p)^{*}=(q,p,p)$ for any projection $(q,p,p)\in {\cal
L}((a,p))$, the $C^{*}$-category $\mathrm{P}(A)$ has enough
isometries, i.e. $\mathrm{P}(A)$ is a pseudo-abelian category.

Let $A$\ be a $C^{*}$-algebra with unit. Denote by $\mathrm{F}(A)$
the additive $C^{*}$-category which has standard Hilbert right
$A$-modules $A^n=A\oplus _{\mathrm{n-times}}\oplus A$ as the objects
and usual $A$-homomorphisms (which has an adjoint) as the
morphisms. Let $\mathrm{P}(A)$ be the standard pseudo-abelian
$C^{*}$-category of $\mathrm{F}(A)$, i.e. $\mathrm{P}(A)=\mathrm{P(F}(A))$.

We say that a subcategory $A$ is a cofinal sub-category of a
$C^{*}$-category $A'$ if for any object $a'\in A'$ there
exists an object $a\in A$ and an isometry $s:a'\rightarrow a$. If
$A$ is an additive full cofinal sub-$C^{*}$-category in a
pseudo-abelian $C^{*}$-category $A'$ then we say that $A'$ {\it is
generated} by the $C^{*}$-category $A$.

\subsection{The $C^*$-categories $Rep(A;B)$ and $\mathrm{Rep}(A;B)$}

We list some examples of $C^*$-categories and $C^*$-categoroids, useful
in the next sections.

1) (See \cite{kas1}) Firstly we define the $C^*$-category
$\mathcal{H}_G(B)$ over a fixed compact second countable group $G$.
The objects of this category are countable generated right Hilbert
$B$-modules equipped with a $B$-linear, norm-continuous $G$-action such
that $g(xb)=g(x)g(b)$ and $<g(x),g(y)>=g<x,y>$, for all $g\in G$.
A morphism $f:E\rightarrow E'$ is $B$-homomorphism commuting
with the action of $G$, such that there
exists $f^*:E'\rightarrow E$ satisfying the conditions:
$<T(x),y>=<x,T^*(y)>$ where $x\in E$ and $y\in E'$. The norm of a
morphism is defined as the norm of linear bounded map. $\mathcal{H}_G(B)$
is an additive $C^*$-category with respect to the sum of the Hilbert modules.
Note that compact group acts on the morphisms by following rule: if
$f:E\rightarrow E'$ then morphism $gf:E\rightarrow E'$ is defined
by formula $gf(x)=g(f(g^{-1}(x)))$ (this action generally
isn't norm-continuous). A morphism is called $invariant$ if $gf=f$.
The category $\mathcal{H}_G(B)$ contains natural class of morphism,
so called compact $B$-homomorphisms \cite{kas1}. Denote it by
$\mathcal{K}_{G}(B)$. Properties of compact $B$-homomorphisms implies
that $\mathcal{K}_{G}(B)$ is a $C^*$-ideal in $\mathcal{H}_G(B)$.
Note that there exists a
$*$-functor $\infty :\mathcal{H}_G(B)\rightarrow \mathcal{H}_G(B)$ and a
natural isomorphism of functors
$id_{\mathcal{H}(B)}\oplus \infty  \simeq \infty $, where
$E^\infty =E\oplus E\oplus \cdots $.
This structure will be called an $\infty $-structure.

2) Now, we define the $C^{*}$-category $\mathrm{rep}_G(A,B)$. The
objects of this category are pairs of form $(E,\phi)$, where $E$ is an
object in $\mathcal{H}_G$ and $\phi :A\rightarrow \mathcal{L}(E)$
is an equivariant $*$-homomorphism. Objects of this type are said to be
$A,B$-bimodules. A morphism $f:(E,\phi )\rightarrow (E',\phi ' )$ is an
$B$-homomorphism $f:E\rightarrow E' $ in $\mathcal{H}_G(B)$ such
that $f\phi (a)=\phi '(a)f$ for all $a\in A$. The structure of
a $C^{*}$-category and action of $G$ is inherited
from the $C^{*}$-category structure of $\mathcal{H}_G(B)$. Let
$\mathrm{rep}(A,B)$ be the sub-$C^*$-category of $\mathrm{rep}_G(A,B)$ which
has the same object as latter but morphisms are invariant under action of
$G$. it is easy to show that $\mathrm{rep}(A,B)$ is an additive
$C^{*}$-category (in fact, a pseudo-abelian $C^{*}$-category). The following
property of $\mathrm{rep}(A,B)$ useful for calculation of the $K$-groups of
$\mathrm{rep}(A,B)$.
The $\infty $-structure on $\mathcal{H}_G(B)$ induces a corresponding
structure on $\mathrm{rep}(A,B)$ via the formulas
$(E,\phi)^\infty =(E^\infty ,\phi ^\infty )$, where
$\phi ^\infty (a)=(\phi (a))^\infty $ for all $a\in A$.

3) Consider the additive $C^{*}$-category $\mathcal{Q}_G(B)$ which is the
quotient $C^{*}$-category $\mathcal{H}_G(B)/\mathcal{K}_G (B)$.
It has an essential compact group action inherited from the action of
a compact group on $\mathcal{H}_G(B)$. Denote by
$\pi :\mathcal{H}_G(B)\rightarrow \mathcal{Q}_G(B)$ the
canonical additive $*$-functor. We need also following
$C^{*}$-category denoted by $\mathcal{Q}_G(A,B)$. By definition objects of
this category have the form $(E,\psi)$, where $E$ is an object in
$\mathcal{H}_G(B)$ and $\psi :A\rightarrow \hom _{\mathcal{Q}_G(B)}(E,E)$
is a equivariantly liftable $*$-homomorphism, i.e., there exists an
$A,B$-bimodule $(E,\phi )$ such that $\pi \phi =\psi $. A morphism
$f:(E,\psi )\rightarrow (E' ,\psi ' )$ is a morphism
$f:E\rightarrow E'$ of the category $\mathcal{Q}_G(B)$ such that
$f\psi(a)=\psi '(a)f$ for all $a\in A$. Let $\mathcal{Q}(B)$ be
sub-$C^*$-category of $\mathcal{Q}_G(B)$ the invariant liftable morphisms
of latter as the morphisms. There is a $*$-functor
$\Theta :\mathrm{rep}(A,B)\rightarrow \mathcal{Q}(A,B)$ given by
$(E,\phi)\mapsto (E', \phi)$ and $f\mapsto \pi (f)$.

4) Now, we want to define the additive $C^{*}$-category
$Rep_G(A,B)$. The class of objects of this category coincides with the class
of objects of $\mathrm{rep}_G(A,B)$. But a morphism
$f:(E,\phi)\rightarrow (E',\phi ')$ is a morphism
$f:E\rightarrow E'$ in $\mathcal{H}(B)$ such that
$$
f\phi (a)-\phi '(a)f\in \mathcal{K}_G(E,E' )
$$
for all $a\in A$. The structure of $C^{*}$-category and action of $G$ are
inherited from $\mathcal{H}_G(B)$.
It is easy to show that $Rep_G(A,B)$ is an additive $C^{*}$-category (but it
isn't a pseudo-abelian $C^{*}$-category). Let $Rep(A;B)$ be
sub-$C^*$-category of $Rep_G(A;B)$ with same class of objects as latter
but morphism are those are invariant under the action of $G$. There is a
canonical additive $*$-functor
$\Pi_{A,B}:Rep(A,B)\rightarrow \mathcal{Q}(A,B)$
defined by
$(E,\phi)\mapsto (E,\pi \phi )$ and $f\mapsto \pi f$.
From the definition follows that the canonical linear map
$$
\hom ((E,\phi ),(E',\phi ' ))\mapsto \hom ((E', \phi ),(E', \phi ' )
$$
is surjective, i.e., $\Pi $ is a Serre functor (see for the
definition \cite{kar}).

There is an $C^*$-ideal $D(A,J;B)$ in $Rep(A;B)$, which is associated
from any closed ideal $J$ in $A$. $C^{*}$-ideal $D(A,J;B)$, which we'll
use in the sequel, is defined by the following manner.
Let $(E,\phi )$ and $(E',\phi ')$ are objects in
$Rep(A,B)$. A morphism $\alpha :(E,\phi )\rightarrow (E',\phi ')$
in $Rep(A,B)$ is in $D(A,J;B)$ if
$\alpha \phi (j)\in \mathcal{K}((E,\phi ),(E'\phi '))$ and
$\phi (j)\alpha \in \mathcal{K}((E,\phi ))$.
The space of all morphisms from $(E,\phi )$ to
$(E',\phi ')$ will be denoted by $D_{\phi ,\phi '}(A,J;E,E';B)$
(if $(E',\phi ')=(E,\phi )$ then it is denoted by $D_\phi (A,J;E;B)$).
Sometimes $Rep(A;B)$ is denoted by $D(A,0;B)$ or $D(A;B)$.

Now we come to our main $C^{*}$-category, that is, $\mathrm{Rep}(A,B)$.

\begin{definition}
\label{q}
Let $\mathrm{Rep}(A,B)$ be the universal pseudo-abelian
$C^{*}$-category of $Rep(A,B)$. Using
the definition of a pseudo-abelian $C^{*}$-category, we have the
following description of $\mathrm{Rep}(A,B)$. Objects of it are triples
$(E,\phi,p)$, where $(E,\phi)$ is an object and
$p:(E,\phi)\rightarrow (E,\phi )$ is a morphism in $Rep(A,B)$ such that
$p^{*}=p$ and $p^2=p$. A morphism
$f:(E,\phi,p)\rightarrow (E' ,\phi ' ,p' )$ is a morphism
$f:(E,\phi )\rightarrow (E,\phi )$ in $Rep(A,B)$ such
that $fp=p'f=f$. In detail, $f$ has the properties
\begin{equation}
f\phi (a)-\phi '(a)f\in \mathcal{K}(E,F),\quad fp=p'f=f.
\end{equation}
\end{definition}

The structure of $C^{*}$-category of $\mathrm{Rep}(A,B)$ comes
from the corresponding structure of $Rep(A,B)$.

\section{Higson's Homotopy Invariance Theorem \label{homin}}

In the next sections we'll need Higson's homotopy invariance
theorem. The purpose of this section is to give account of Higson's
homotopy invariant theorem, both for real and complex cases.
Higson's theorem asserts the following. Let $E$ be a stable and
split additive functor from admissible sub-category of the category $complex$
$C^{*}$-algebras into the category of abelian groups. Then it is homotopy
invariance \cite{hig3} (Theorem 3.2.2).

This important theorem plays major role for setting homotopy invariance of
functors (for example, using this theorem, it was solved well known problem
about isomorphism of algebraic and topological $K$-theories on stable
$C^{*}$-algebras \cite{suw}). The natural question arises whether this
theorem is for real $C^{*}$-algebras? The reason for it is that we couldn't
dissemination the proof of Lemma 3.1.2 and Theorem 3.11 of \cite{hig3} for
real $C^{*}$-algebras; because the complexity of $C^{*}$-algebras is
significant to prove them. We choose another point of view for Higson's
theorem for both cases mentioned above, namely, an investigation of
$KK$-theory by J. Cuntz and G. Scandalis \cite{cusk}. We have chosen this
interpretation first, because it is simple and only very simple arguments
are used from $KK$-theory, and on the other hand it is right both for
complex and real $C^{*}$-algebras. Equivalence relation on Kasparov
bimodules is defined in this article that is called ''Cobordism''. It's main
feature is that this equivalence coincides with the Kasparov's induced by
homotopy of bimodules.

Let $(E,\varphi ,F)$ be a Kasparov $A,B$-bimodule, where $E$ is a countable
generated Hilbert right $B$-module, $F\in \mathcal{L}_B(E)$ is a degree one
element and $\varphi :A\rightarrow \mathcal{L}_B(E)$ is a
$*$-homomorphism.  We
say that $A,B$-bimodules $(E_0,\varphi _0,F_0)$ and $(E_1,\varphi _1,F_1)$
are $cs$-isomorphic if there exists such a degree zero unitary
$u:E_0\rightarrow E_1$ that $u\varphi _0(a)u^{*}=\varphi _1(a)$
for any $a\in A$ and $uF_0u^{*}-F_1\in \mathcal{K}(E_1)$. Now, one
easily checks that two Kasparov $A,B$-bimodules $(E_0,\varphi _0,F_0)$
and $(E_1,\varphi _1,F_1)$ are cobordant if and only if there exists such
a Kasparov $A,B$-bimodule
$(E,\varphi ,F)$ that $(E_0,\varphi _0,F_0)\oplus (E,\varphi ,F)$ and
$(E_1,\varphi _1,F_1)\oplus (E,\varphi ,F)$ are $cs$-isomorphic. Thus from
\cite{cusk} (Theorem 3.7) one deduces as follows (cf. remark 3.8 in
\cite{cusk}).

{\it Let} $A$ {\it be a separable} $Z_2$-{\it graded} $C^{*}$-{\it algebra
and} $B$ {\it be a }$Z_2${\it -graded } $\sigma ${\it -unital }
$C^{*}${\it -algebra. Let} $\overline{KK}(A,B)$ {\it be the cancellation
monoid of the abelian monoid of classes} $A,B${\it -bimodules identified
by} $cs${\it -isomorphism and }$KK(A,B)$ {\it be Kasparov group.
Then the natural homomorphism}
$$
\tau :\overline{KK}(A,B)\rightarrow KK(A,B)
$$
{\it is an isomorphism. In particular}, $\overline{KK}(A,B)$
{\it is an abelian group and homotopy invariance for both arguments}.

Bellow, we give a slightly different variant of this theorem, which will be
our fulcrum to approach Higson's theorem.

Consider Kasparov $A,B$-bimodules satisfying conditions $F^{*}=F$ and
$F^2=1$ (we call them $special$). Let $\widehat{KK}(A,B)$ be the
cancellation monoid of the abelian monoid of classes special $A,B$-bimodules
identified by $cs$-isomorphism. Then one formulate the above mentioned result in the following
form.

$\mathbf{Theorem}$ A. $Let$ $A$ $be$ $a$ $separable$ $Z_2$-$graded$
$C^{*}$-$algebra$ $and$ $B$ $it$ $be$ $a$ $Z_2$-$graded$ $\sigma $
 -$unital$ $C^{*}$
 {\it -algebra. Let }$\widehat{KK}(A,B)${\it be the
cancellation monoid of the abelian monoid of classes special}
$A,B${\it -bimodules identified by} $cs${\it -isomorphism and }$KK(A,B)$
{\it be Kasparov group. Then the natural homomorphism}
$\tau :\widehat{KK}(A,B)\rightarrow KK(A,B)$ {\it \ is an isomorphism.
In particular}, $\widehat{KK}(A,B)$ {\it is an abelian group and homotopy
invariance for both arguments}.

The proof of this theorem is based on the fact that cobordism and homotopy
of Kasparov bimodules coincides, and from Lemmas 17.4.2-17.4.3 in \cite{bla}.
We shortly remind the content of them:

Kasparov $KK$-groups won't be changed if in their definition will be taken
only

(1) Kasparov bimodules with property $F^{*}=F$, and one takes only homotopy
or only cobordism;

(2) Kasparov bimodules with property $F^{*}=F$ and $F^2=1$, and one takes
homotopy or only cobordism;

(3) special Kasparov bimodules, and one takes only $cs$-isomorphism.

We remark that the functional calculus of a self-adjoint element is used to
show (1)-(3) in \cite{bla}. The functional calculus for a self-adjoint
element $x$ in complex $C^{*}$-algebra $A$ is the $*$-monomorphism
$\Phi :C(\mathrm{sp}x)\rightarrow A$, defined by the rule
$id_{\mathrm{sp}x}\mapsto x$. For the case of a real $C^{*}$-algebra under
functional calculus we mean the following. Let $A$ be a real $C^{*}$-algebra, and consider the complex
involutive algebra $A\otimes _RC$ with involution
$(a\otimes c)^{*}=a^{*}\otimes \bar c$. Then there exists a $C^{*}$-norm,
with respect to which $A\otimes _RC$ is a complex $C^{*}$-algebra, and
canonical $*$-embedding $i:A\rightarrow A\otimes _RC$ defined by
$a\mapsto a\otimes 1$ is an isometry (cf. Theorem 2 and Corollary 2
in \cite{pal}, also \cite{kas1}).
For real $C^{*}$-algebras there exists functional calculus for
self-adjoint elements. To be more precise, let $r\in A$ is a self-adjoint
element and $A(r)$ be the real closed sub-algebra in $A$ generated by $r$.
It is real part of the sub-$C^{*}$-algebra $C^{*}(r\otimes 1)$, generated
by the element $r\otimes 1$ in $A\otimes _RC$. According to functional
calculus of complex $C^{*}$-algebras one has the $*$-monomorphism
$A(r)\rightarrow A\otimes _RC$ defined by rule $r\rightarrow r\otimes 1$.
Thus under functional calculus we mean the functional calculus of then
element $r\otimes 1$ in $A\otimes _RC$ (For example, if $f$ is
continuous real function on the $\mathrm{sp}(r\otimes 1)$ then there
exists the unique element $f(r)\in A(r)$ that $i(f(r))=f(r\otimes 1)$).

Now theorem A may be interpreted for trivially graded $A$ and $B$ in the
following way.

Let $(\varphi ,\psi ,U)$ be a triple, where $H_\varphi $ and $H_\psi $ are
countable generated Hilbert right $B$-modules with trivial grading; $\varphi
:A\rightarrow \mathcal{L}_B(H_\varphi )$\ and
$\psi :A\rightarrow \mathcal{L}_B(H_\psi )$\ are $*$-homomorphisms;
$U:H_\varphi \rightarrow H_\psi $ is a $B$-homomorphism which has an adjoint.
A triple $(\varphi ,\psi ,U)$ is said to be {\it unitary} {\it Fredholm}
$A,B${\it -bimodule} if $U$ is a unitary and satisfies the following
condition:

$U\varphi (a)-\psi (a)U\in \mathcal{K}(H_\varphi ,H_\psi )$, for all
$a\in A$.

$A,B$-bimodules $(\varphi _0,\psi _0,U_0)$ and
$(\varphi _1,\psi _1,U_1)$ are said to be  $cs$-isomorphic if there
exists degree zero unitaries
$u:H_{\varphi _0}\rightarrow H_{\varphi _1}$ and
$v:H_{\psi _0}\rightarrow H_{\psi _1}$ that
$u\varphi _0(a)u^{*}=\varphi _1(a)$,
$v\psi _0(a)v^{*}=\psi _1(a)$
and
$$
vU_0u^{*}-U_1\in \mathcal{K}(H_{\varphi _1},H_{\psi _1})
$$
for any $a\in A$.

The sum of $A,B$-bimodules are defined in usual manner
$$
(\varphi _0,\psi
_0,U_0)\oplus (\varphi _1,\psi _1,U_1)=(\varphi _0\oplus \varphi _1,\psi
_0\oplus \psi _1,U_0\oplus U_1).
$$

In the definition of group $\widehat{KK}(A,B)$ Kasparov special
$A,B$-bimodules can be replaced by Fredholm unitary $A,B$-bimodules. This
statement is proved by the following way. First of all, if $B$ is trivially
graded $C^{*}$-algebra then any graded Hilbert $B$-module $E$ can be
represented in this form $E=E_0\oplus \bar E_1$ where $E_0$ and $E_1$ are
trivially graded Hilbert $B$-modules and $\bar E_1$ is the opposite to $E_1$
graded Hilbert $B$-module, i.e. $\bar E_1^{(0)}=0$ and $\bar E_1^{(1)}=E_1$.
Let $i_d:E_1\rightarrow \bar E_1$ be the degree one $B$-linear map defined
by the identity map. Now, if $(E,\phi ,F)$ is a Kasparov special $A,B$%
-bimodule then one has
\begin{equation}
\label{kafre}E=E_0\oplus \bar E_1,\;\;\;\phi =\left(
\begin{array}{cc}
\varphi & 0 \\
0 & \bar \psi
\end{array}
\right) ,\;\;\;{\rm and}\;\;\;F=\left(
\begin{array}{cc}
0 & \overline{U}^{*} \\ \overline{U} & 0
\end{array}
\right) .
\end{equation}
Here $\varphi :A\rightarrow \mathcal{L}_B(E_0)$ and
$\psi :A\rightarrow \mathcal{L}_B(E_1)$ are $*$-homomorphisms,
$U:E_0\rightarrow E_1$ is unitary $B$-homomorphism,
$\overline{U}=i_d\cdot U\cdot i_d^{*}$ and $\bar \psi (a)=i_d\cdot \psi (a)
\cdot i_d^{*}$. It is easily check that $(\varphi ,\psi ,U)$ is
Fredholm unitary $A,B$-bimodule. Conversely, let $(\varphi ,\psi ,U)$
be a Fredholm unitary $A,B$-bimodule, where
$\varphi :A\rightarrow {\cal L}_B(E_0)$
and $\psi :A\rightarrow {\cal L}_B(E_1)$ are $*$-homomorphisms such that
$E_0$ and $E_1$ are trivially graded Hilbert $B$-modules. Then the
correspondence Kasparov special $A,B$-bimodule $(E,\phi ,F)$ is given by the
formulas \ref{kafre}. Therefore Theorem A may be formulated as follows.

\begin{theorem}
\label{csf} Let $A$\ be a separable trivially graded $C^{*}$-algebra and
$B$ be a trivially graded $\sigma $-unital $C^{*}$-algebra. Let
$\widetilde{KK}(A,B)$ be the cancellation monoid of the abelian monoid of
classes Fredholm $A,B$-bimodules identified by $cs$-isomorphism, and
$KK(A,B)$\ be Kasparov group. Then the natural homomorphism
$$
\tilde \tau :\widetilde{KK}(A,B)\rightarrow KK(A,B)
$$
is an isomorphism. In particular,$\widetilde{KK}(A,B)$\ is an abelian
group and homotopy invariance for both arguments.
\end{theorem}

Let $\mathcal{S}_0$ be be a sub-category of $C^{*}$-algebras containing
two $C^{*}$-algebras, $\mathbf{F}$ and $\mathbf{F}[0;1]$, as the objects;
and two non-trivial morphisms, these are evolution maps
$ev_0:\mathbf{F}[0;1]\rightarrow \mathbf{F}$ and
$ev_1:\mathbf{F}[0;1]\rightarrow \mathbf{F}$ given
by the formulas $ev_0(f(t))=f(0)$ and $ev(f(t))=f(1)$, where
$f(t)\in \mathbf{F}[0;1]$. The following definition points out an important
class of functors on $\mathcal{S}_0$.

\begin{definition}
\label{wh} A functor $E$ from $\mathcal{S}_0$ into a category is said to be
\\
(1) {\it weak homotopy} (or shortly $w${\it -homotopy})
if $E(ev_0)=E(ev_1)$.
\\
(2) {\em weak} $K$ {\em -hereditary} if there is
a natural transformation
of functors
$$
\chi :\widetilde{KK}(-;\mathbf{F})\rightarrow \hom (E(-);E(\mathbf{F}))
$$
so that $\chi (e,\theta ,1)$ is the identity morphism from $E(\mathbf{F})$
in itself, where $(e,\theta ,1)$ is Fredholm unitary $\mathbf{F,F}$-module;
$e:\mathbf{F}\rightarrow \mathcal{K}(\mathcal{H})$ is a $*$-homomorphism
which maps $1\in \mathbf{F}$ in a rank one projection,
$\theta :\mathbf{F}\rightarrow \mathcal{K}(\mathcal{H}$ is the zero
homomorphism and $1$ is unit of the algebra $\mathcal{L}(\mathcal{H})$,
where $\mathcal{H}$ is  a countable generated Hilbert space.
\end{definition}

For example, let $E$ be a homotopy invariant covariant functor
on the category $S$ of separable $C^*$-algebras and $*$-homomorphisms.
Then functor $E(A\otimes -)$ is a $w$-homotopy functor on $S_0$,
for any $A\in S$. Conversely, if $E(A\otimes -)$ is a $w$-homotopy functor
on $S_0$, for any $A\in S$ then, of course, $E$ is homotopy invariant too.

Next we'll only need the following very particular
case of Theorem \ref{csf}.

\begin{corollary}
\label{pccsf}The contravariant functor
$\widetilde{KK}(-;\mathbf{F})$ on the category $\mathcal{S}_0$ is
$w$-homotopy invariant.
\end{corollary}

Now, from Corollary \ref{pccsf} we deduce the following lemma.
We'll use this lemma to obtain Higson's theorem.

\begin{lemma}
\label{hw}Let $E$ be a weak $K$-hereditary functor from $S_0$ into a
category. Then $E$ is $w$-homotopy functor.
\end{lemma}

\begin{proof} Since $K$-homology, by Corollary \ref{pccsf}, is $w$-homotopy
and $E$ is a weak $K$-hereditary on $S_0$, the diagram
$$
\begin{array}{ccc}
\widetilde{KK}(\mathbf{F,F}) & \stackrel{\chi _{\mathbf{F}}}{\rightarrow } &
\hom (E(\mathbf{F}),E(\mathbf{F})) \\ \;\;\;\;\;\;\downarrow ^{e\tilde v_0=e
\tilde v_1} & & \;\;\downarrow ^{e\bar v_0}\downarrow ^{e\bar v_1} \\
\widetilde{KK}(\mathbf{F}[0;1],\mathbf{F}) &
\stackrel{\chi _{\mathbf{F}[0,1]}}{\rightarrow } &
\hom (E(\mathbf{F}[0;1]),E(\mathbf{F}))
\end{array}
$$
commutes. Thus $e\bar v_0\cdot \chi _{\mathbf{F}}=e\bar v_1\cdot
\chi _{\mathbf{F}}$. Again, since $E$ is a weak $K$-hereditary, one has
$E(ev_0)=e\bar v_0(\chi _{\mathbf{F}}(\iota ))=e\bar
v_1(\chi _{\mathbf{F}}(\iota ))=E(ev_1)$
where $\iota $ is the class of Fredholm $\mathbf{F,F}$-bimodule
$(e,\theta ,1)$. Therefore $E$ is a $w$-homotopy functor.
\end{proof}

\subsection{Pairings with Fredholm Pairs}

Let $E:S\rightarrow Ab$ be a functor, where $Ab$ is the category abelian groups
and homomorphisms. A pairing of $E$ with the set of Fredholm pairs is defined
in \cite{hig3}. We'll recall it.

A Fredholm $B$-pair is such a pair $(\varphi ,\psi )$ where $\varphi $ and
$\psi $ are $*$-homomorphisms from $B$ into
$\mathcal{L}_{\mathbf{F}}(\mathcal{H})$
that $\varphi (b)-\psi (b)\in \mathcal{K}(\mathcal{H})$ for any $b\in B$,
where $\mathcal{H}$ is a countable generated Hilbert space over
$\mathbf{F}$, Here $\mathcal{K}(\mathcal{H})$ is the $C^{*}$-algebra of
compact operators. A pairing of $E$ with the set of Fredholm $B$-pairs is
a rule. This assigns to each Fredholm $B$-pair $(\varphi ,\psi )$ a morphism
$\times (\varphi ,\psi ):E(A\otimes B)\rightarrow E(A\otimes \mathbf{F})$
in the category $C$, for any $A\in S$ and $B\in S_0$ , with the following
properties:

\begin{enumerate}
\item  \underline{Functoriality. } If $(\varphi ,\psi )$\ is a Fredholm
$B'$-module, and if $f:B\rightarrow B'$ is a $*$-homomorphism from
$S$, then the diagram
$$
\begin{array}{ccc}
E(A\otimes B) & \stackrel{\times (\varphi f,\psi f)}{\longrightarrow}&
E(A\otimes \mathbf{F})\\
\downarrow ^{E(\mathrm{id}_A\otimes f)}& &\parallel \\
E(A\otimes B') & \stackrel{\times (\varphi ,\psi )}{\longrightarrow}&
E(A\otimes \mathbf{F})
\end{array}
$$
commutes.

\item  \underline{Additivity.} If $(\varphi ,\chi )$\ and $(\chi ,\psi )$\
are Fredholm $B$-pairs, then
$$
\times (\varphi ,\chi )+\times (\chi ,\psi )=\times (\varphi ,\psi ).
$$

\item  \underline{Stability.} If $(\varphi ,\psi )$\ is a Fredholm
$B$-pair and $\eta :B\rightarrow \mathcal{L}(\mathcal{H})$ is any
$*$-homomorphism then
$$
\times (\varphi ,\psi )=\times \left( \left(
\begin{array}{cc}
\varphi  & 0 \\
0 & \eta
\end{array}
\right) ,\left(
\begin{array}{cc}
\psi  & 0 \\
0 & \eta
\end{array}
\right) \right) .
$$

\item  \underline{Non-degeneracy.} If $(e,\theta)$\ is a Fredholm
$\mathbf{F}$-module for which
$e:\mathbf{F}\rightarrow \mathcal{K}(\mathcal{H})$
maps $1\in \mathbf{F}$ to $p$, where $p$ is a rank
one projection in $\mathcal{K}(\mathcal{H})$ and $\theta $ is the zero
homomorphism. Then
$$
\times (e,0):E(A\otimes \mathbf{F})\rightarrow E(A \otimes \mathbf{F})
$$
is the identity morphism.

\item  \underline{Unitary equivalence.} If $U\in \mathcal{L}(\mathcal{H})$
is a unitary then
$$
\times (\varphi ,\psi )=\times (U\varphi U^{*},U\psi U^{*})
$$

\item  \underline{Compact perturbations.} If $U\in \mathcal{L}(\mathcal{H})$
is a unitary which is equal to the identity modulo compacts, then
$$
\times (\varphi ,U\varphi U^{*})=0.
$$
\end{enumerate}

The following theorem is exactly theorem 3.1.4 of \cite{hig3} (not only on
the category of complex $C^*$-algebras as in \cite{hig3}, but on
the category of real $C^*$-algebras too.)

\begin{theorem}
\label{frp}
Let be $E$ a functor from the category $\mathcal{S}$ into the category $Ab$
of abelian groups and their homomorphisms and the functor admits a pairing
with the set of Fredholm $B$-pairs, $B\in S_0$. Then $E$ is a homotopy functor.
\end{theorem}

\begin{proof} As it was pointed out after definition \ref{wh}, it is
enough to show that $E(A\otimes -)$ is a weak $K$-hereditary functor.
So, we'll have to construct a natural transformation functors
$$
\chi :\widetilde{KK}(-;\mathbf{F})\rightarrow \hom
(E(A\otimes -);E(\mathbf{A\otimes F})).
$$
Let $(\varphi ,\psi ,U)$ be a Fredholm unitary $B,\mathbf{F}$-bimodule,
then by definition
$$
\chi (\varphi ,\psi ,U)=\times (\varphi ,U^{*}\psi U).
$$
The correspondence $\chi $ is correctly defined, since
$\widetilde{KK}(-;\mathbf{F})$ is the cancellation monoid of the monoid of
the classes of $cs$-isomorphic unitary $A,\mathbf{F}$-bimodules and the
following hold.

(1). {\it If }$(\varphi _0,\psi _0,U_0)${\it and }
$(\varphi _1,\psi_1,U_1) $ {\it are }$cs${\it -isomorphic then}
$\chi (\varphi _0,\psi _0,U_0)=\chi (\varphi _1,\psi _1,U_1)$.

Indeed, let $(u,v)$ be a $cs$-isomorphism from
$(\varphi _0,\psi _0,U_0)$ in $(\varphi _1,\psi _1,U_1)$, then
$$
\begin{array}{c}
\chi (\varphi _0,\psi _0,U_0)=\times (\varphi _0,U_0^{*}\psi _0U_0)=(
\mathrm{by\;unitary\;equivalence})\\= \times (\varphi _1,u^{*}(U_0^{*}\psi
_0U_0)u)=\times (\varphi _1,(u^{*}U_0^{*}v)(v^{*}\psi _0v)(v^{*}U_0u))=(
\mathrm{by\;additivity}) \\ =\times (\varphi _1,U_1^{*}\psi _1U_1)+\times
(U_1^{*}\psi _1U_1,(u^{*}U_0^{*}v)\psi _1(v^{*}U_0u)))=(
\mathrm{by\;compact\;perturbation}) \\ =\times (\varphi _1,U_1^{*}\psi
_1U_1)=\chi (\varphi _1,\psi _1,U_1).
\end{array}
$$

(2). $\chi (\varphi _0\oplus \varphi _1,\psi _0\oplus \psi _1,U_0\oplus
U_1)=\chi (\varphi _0,\psi _0,U_0)+\chi (\varphi _1,\psi _1,U_1)$.

Indeed,
$$
\begin{array}{c}
\chi (\varphi _0\oplus \varphi _1,\psi _0\oplus \psi _1,U_0\oplus U_1)=(
\mathrm{by\;additivity}) \\ =\chi (\varphi _0\oplus \varphi _1,\psi _0\oplus
\varphi _1,U_0\oplus 1)+\chi (U_0^{*}\psi _0U_0\oplus \varphi _1,U_0^{*}\psi
_0U_0\oplus \psi _1,1\oplus U_1)=(\mathrm{by\;stability})
\\ =\chi (\varphi _0,\psi _0,U_0)+\chi (\varphi _1,\psi _1,U_1)
\end{array}
$$
\end{proof}

Now, Higson's theorem may be formulated in the following way.

\begin{theorem}
\label{higt}
Let $E$ be stable and split exact covariant or contravariant functor on
the category of separable complex or real $C^*$-algebras. Then $E$ is
homotopy invariant.
\end{theorem}

\begin{proof} The complex case, using theorem 3.1.4, was proved in
\cite{hig3}. The latter may be applied, mutatis mutandis, to the real
case, according to Theorem \ref{frp}. The contravariant case can be deduced
from the covariant case in following manner. Let $\Phi $ be a
contravariant functor with above mentioned properties and $A$ be a separable
$C^{*}$-algebra, then
consider a covariant functor $\hom (\Phi (-),\Phi (A))$ which, of
course, is split exact and stable and thus homotopy invariant.
Therefore,
$$
\Phi (ev(0))=\Phi (ev(0))(id_{\Phi (A))}=\Phi (ev(1))(id_{\Phi
(A))}=\Phi (ev(1)),
$$
where $ev(0),ev(1):A([0;1])\rightarrow A$ are evolutions at $0$
and $1$.
\end{proof}

Let $S_G$ be the category of separable $C^*$-algebras with action of
a fixed compact second countable group as the objects and the equivariant
$*$-homomorphisms as the morphisms. For a fixed $B\in S_G$, one has natural
functor $B\otimes -:S\rightarrow S_G$ sending $C^*$-algebra A in the
$B\otimes A$ with diagonal action of $G$. Trivial checking shows that latter
functor sends a split exact sequence from $S$ in the split exact sequence of
the category $S_G$.

Further, consider algebra $\mathcal{K}$ as the object in $S_G$, assuming that
action of $G$ on $\mathcal{K}$ is trivial. We say that functor $E$ is stable
if $E(e_A)$ is an isomorphism, where $e_A$ is a natural equivariant inclusion
$e_A:A\rightarrow A\otimes \mathcal{K}$ given by the map
$a\mapsto a\otimes p$, $p$ is rank one projection.

Now, we'll make some remarks on the notation of homotopy. Let $A$ be an algebra
in $S_G$. Consider algebra $A[0;1]$ with action of $G$ defined by equality
$(g\cdot f)(t)=g(f(t)), \forall t\in [0;1]$. It is well known that the latter
algebra is isomorphic to the tensor product $A\otimes \mathbf{F}[0;1]$ where
the action of $G$ on $\mathbf{F}[0;1]$ is trivial. Let $E$ be a functor on
$S_{G}$. It is easily to show that $E$ is homotopy invariant if only if
$E(e_A(0))=E(e_A(1))$ for any $A\in S_G$, where $e_t:A[0;1]\rightarrow A$
the evolution at $t\in [0;1]$. This fact may be
trivially reformulated in the following way. A functor $E$ on $S_G$ is
homotopy invariant if and only if the functor $E(A\otimes -)$ is homotopy
invariant on the category $S$ for all $A\in S_{G}$.

A functor $E$ will be said {\it split exact} if a functor  $E(A\otimes -)$ is
split exact on $S$, for any algebra $A$ in $S_G$.

The following corollary is a consequence of Theorem \ref{higt}.

\begin{corollary}
\label{ghigt}
Let $E$ be stable and split exact functor on the category $S_G$ of real or
complex $C^*$-algebras with the action of a fixed compact second countable
group $G$ and equivariant $*$-homomorphisms. Then $E$ is homotopy invariant.
\end{corollary}

\begin{proof} For any object $A$ in $S_G$, consider the functor
$E(A\otimes -)$ which, of course is stable and split exact on the category
$S$, thus homotopy invariant by Theorem \ref{higt}. Therefore $E$ is
homotopy invariant on the category $S_G$, by the principle pointed after
Theorem \ref{higt}.
\end{proof}

 \section{On the Cuntz-Bott Periodicity \label{cbcoho}}

As before, let $S_G$ denotes category of separable complex or real
$C^*$-algebras with action of compact second countable group $G$
and equivariant $*$-homomorphisms. let
$\mathbb{E}=\{E_n\}_n\in mathbf{Z}$ be a set of covariant functors
from $S_G$ into the category of abelian groups and homomorphisms,
indexed by the integer numbers. One says that $\mathbb{E}$ is Cuntz-Bott
homology theory on the category $\mathbb{E}$ if

\begin{enumerate}
  \item $\mathbb{E}$ has weak excision property. Namely, for any exact
  sequence $0\rightarrow I\rightarrow B\rightarrow A\rightarrow 0$ of
  algebras from $S_G$, where epimorphism admits equivariant completely
  positive an contracive section. Then

  i. there exist homomorphism $\delta _n : E_n(A)\rightarrow E_{n-1}(I)$,
  for any $n\in Z$, non-depending on a completely positive
  and contractive section of $p$, and natural in the following sense.
  let
  $$
  \begin{array}{ccccccccc}
    0 & \rightarrow & I& \rightarrow & B & \rightarrow & A& \rightarrow
     & 0 \\
     & & \;\;\; \downarrow ^{f_I}&  & \;\;\;\downarrow ^{f_B}&
     & \;\;\;\downarrow ^{f_A}& &  \\
    0 & \rightarrow& I' & \rightarrow & B' & \rightarrow& A' & \rightarrow
     & 0 \
  \end{array}
  $$
   be a commutative diagram such that in the horizontal short exact
   sequences epimorphisms have completely positive and contractive sections.
   Then diagram
   $$
   \begin{array}{ccc}
     E_n(A) & \stackrel{\delta}{\rightarrow} & E_{n-1}(I) \\
    \; \;\;\downarrow ^{E_n(f_A)}&  & \;\;\;\downarrow ^{E_n(f_I)}\\
     E_n(A') & \stackrel{\delta '}{\rightarrow}& E_{n-1}(I') \\
   \end{array}
   $$
   commutes.
   ii. The natural two-sided sequence of abelian groups
   $$
   \cdots \rightarrow E_n(I)\rightarrow E_n(B)\rightarrow
   E_n(A)\stackrel{\delta _n}{\rightarrow }E_{n-1}(I)\rightarrow \cdots
   $$
   is exact.
  \item $\mathbb{E}$ is stable. This means that if
  $e_A:A\rightarrow A\otimes \mathcal{K}$ is a homomorphism defined by a map
  $a\mapsto a\otimes p$, here $p$ is a rank one projection in $\mathcal{K}$,
  then $E_n(e_A)$ is an isomorphism.
\end{enumerate}
Now, we make some remarks about Cuntz's results on Bott periodicity.
Let $\mathbf{T}_{\mathbb{C}}$ Toeplitz complex $C^*$-algebra generated
by an isometry. Denote by $C_{\mathbb{C}}(S^1)$ the $C^*$-algebra of
continuous complex functions on the standard unite cycle $S^1$ of
the module one complex numbers. There is a short exact sequence
$$
0\rightarrow \mathcal{K}_{\mathbb{C}}\rightarrow \mathbf{T}_{\mathbb{C}}
\stackrel{t}{\rightarrow }C_{\mathbb{C}}(S^1)\rightarrow 0.
$$

The real case may be considered in the following way. There is a "Real"
structure on $\mathbf{T}_{\mathbb{C}}$ defined by equality $\bar{v}=v$.
Similarly, on the $C_{\mathbb{C}}(S^1)$ a "Real" structure is defined by the map
$f(z)\mapsto \overline{f(\bar{z})}$. Denote by $C_{\mathbb{R}}(S^1)$ the
real sub-algebra of fixed elements in $C_{\mathbb{C}}(S^1)$ relative the latter
conjugation. Let
$t:\mathbf{T}_{\mathbb{R}}\rightarrow C_{\mathbb{R}(S^1)}$
is given by a map $v\mapsto \mathrm{id}_{S^1}$. One has a short exact sequence
$$
0\rightarrow \mathcal{K}_{\mathbb{R}}\rightarrow \mathbf{T}_{\mathbb{R}}
\stackrel{t}{\rightarrow }C_{\mathbb{R}}(S^1)\rightarrow 0.
$$

Note that natural projection $p:\mathbf{T}_{\mathbb{F}}\rightarrow \mathbb{F}$
defined by the map $v\mapsto 1$ splits by the map
$j:\mathbb{F}\rightarrow \mathbf{T}_{\mathbb{F}}$
defined by $1\mapsto1$. Let $\mathbf{T'}_{\mathbb{F}}$ be the kernel of $p$.
The following proposition is completely analogue of the Proposition 4.3.
\cite{cunt}.
\begin{proposition}
\label{ttob}
Let $\mathbb{E}$ be a Cuntz-Bott homology. Then the homomorphisms
$E_n(\mathrm{id}_A\otimes j)$ and $E_n(\mathrm{id}_A\otimes p)$ are
isomorphisms between $E_n(A)$ and
$E_n(A\otimes \mathbf{T}_{\mathbb{F}})$ for arbitrary $A\in S_G$ and $n\in Z$. In
particular, $E_n(A\otimes \mathbf{T'}_{\mathbb{F}})=0$.
\end{proposition}

\begin{proof}
The sequence
$$
0\rightarrow A\otimes \mathbf{T}'_{F}\rightarrow A\otimes
\mathbf{T}_{\mathbb{F}}\rightarrow A\otimes \mathbb{F}\rightarrow 0
$$
is split exact. Since functors $E_n$ are split exact and stable, they are
homotopy invariant by Higson's theorem. Then the proof literally coincides
with the proof of Proposition 4.3. in \cite{cunt}.
\end{proof}

Now, let $\mho _{\mathbb{F}}$ be a sub-algebra in $C_{\mathbb{F}(S^1)}$
making sequence
$$
0\rightarrow \mho _{\mathbb{F}}\rightarrow C_{\mathbb{F}(S^1)}\rightarrow
\mathbb{F}\rightarrow 0
$$
exact. Then one has natural exact sequence
$$
0\rightarrow \mathcal{K} _{\mathbb{F}}\rightarrow
\mathbf{T}'_{\mathbb{F}}\rightarrow
\mho _{\mathbb{F}}\rightarrow 0.
$$
From the definition of $\mho _{\mathbb{F}}$ it follows that latter algebra
is nuclear $C^*$-algebra. This implies that epimorphism in the latter
exact sequence has a completely positive and contractive section. Thus the
epimorphism has completely positive and contractive section in the exact
sequence
$$
0\rightarrow A\otimes \mathcal{K}_{\mathbb{F}}\rightarrow
A\otimes \mathbf{T}'_{\mathbb{F}}\rightarrow
A\otimes \mho _{\mathbb{F}}\rightarrow 0,
$$
by Lemma 1.3.4. in \cite{hig3}. It implies, According to the Proposition
\ref{ttob},
that $E_n(A)$ is naturally isomorphic to
$E_{n+1}(A\otimes \mho _{\mathbb{F}})$. Summarize this fact we have the following.

\begin{theorem}
\label{ztnuc}
(cf. $\mathrm{Cuntz}$, \cite{cunt}) Let $\mathbb{E}$ be a Cuntz-Bott homology
theory. Then there is a natural isomorphism
\begin{equation}
\label{isocunt}
E_n(A)=E_{n+1}(A\otimes \mho _{\mathbb{F}}).
\end{equation}
\end{theorem}

From this useful theorem we deduce an elementary but applicable principle
which we'll use in the sequel.

\begin{corollary}
\label{natiso}
Let $\mathbb{E}$ and $\mathbb{E}'$ are Cuntz-Bott homology theories. If
there exists natural isomorphism $\mu _m:E_m\rightarrow E'_m$ for some $m\in Z$,
then there exists natural isomorphism $\mu _n:E_n\rightarrow E'_n$ for all
$n\in Z$.
\end{corollary}

\begin{proof}
Consider the following short exact sequence
$$
0\rightarrow \mathbf{F}(0;1)\rightarrow
\mathbf{F}(0;1]\rightarrow \mathbf{F}\rightarrow 0
$$
where $\mathbf{F}$ is $\mathbb{R}$ or $\mathbb{C}$ the fields
of real or complex numbers respectively. Since each algebra in the short
exact sequence are nuclear $C^*$-algebras, according to Lemma 1.3.4.
\cite{hig3}, one concludes that in the exact sequence
$$
0\rightarrow A\otimes \mathbf{F}(0;1)\rightarrow
A\otimes \mathbf{F}(0;1]\rightarrow A\otimes \mathbf{F}\rightarrow 0
$$
the epimorphism has a completely positive and contractive section,
for any separable $C^*$-algebra $A$. From
definition of Cuntz-Bott homology theory and Higson's homotopy invariant
theorem immediately follows that any functor $E_n$ $E'_n$ are homotopy
invariant. From the latter short exact sequence it follows that there
are natural isomorphisms
\begin{equation}
\label{isosusp}
E_{n+1}(A)\simeq  E_n(A\otimes \Omega _{\mathbb{F}})\;\;\;\;\mathrm{and}\;\;\;\;
E'_{n+1}(A)\simeq  E'_n(A\otimes \Omega _{\mathbb{F}}),
\end{equation}
where $\Omega _{\mathbb{F}}=\mathbf{F}(0;1)$. On the one hand, the formulas
\ref{isosusp} guarantees natural isomorphism $E_n(A)\simeq E'_n(A)$ for
$n\geq m$. On the other hand, Cuntz isomorphisms \ref{isocunt} guarantees
natural isomorphism $E_n(A)\simeq E'_n(A)$ for
$n\leq m$.
\end{proof}

The definition of Cuntz-Bott cohomology theory is dual to homology theory
case. Let $\mathbb{E}=\{E^n\}_{n\in Z}$ be a Cuntz-Bott cohomology theory,
We use the following identification $E^n=E_{-n}$ in sequel. Then Definition
of Cuntz-Bott cohomology theory has the following form.

Let, as before, $S_G$ denotes category of separable complex or real
$C^*$-algebras with action of compact second countable group $G$
and equivariant $*$-homomorphisms. let
$\mathbb{E}=\{E_n\}_n\in \mathbf{Z}$ be a set of contravariant functors
from $S_G$ into the category of abelian groups and homomorphisms,
indexed by the integer numbers. One says that $\mathbb{E}$ is Cuntz-Bott
cohomology theory on the category $\mathbb{E}$ if

\begin{enumerate}
  \item $\mathbb{E}$ has weak excision property. Namely, for any exact
  sequence $0\rightarrow I\rightarrow B\rightarrow A\rightarrow 0$ of
  algebras from $S_G$, where epimorphism admits equivariant completely
  positive an contracive section. Then

  i. there exist homomorphism $\delta _n : E_n(I)\rightarrow E_{n-1}(A)$,
  for any $n\in Z$, non-depending on a completely positive
  and contractive section of $p$, and natural by the following sense.
  let
  $$
  \begin{array}{ccccccccc}
    0 & \rightarrow & I& \rightarrow & B & \rightarrow & A& \rightarrow
     & 0 \\
     & & \;\;\; \downarrow ^{f_I}&  & \;\;\;\downarrow ^{f_B}&
     & \;\;\;\downarrow ^{f_A}& &  \\
    0 & \rightarrow& I' & \rightarrow & B' & \rightarrow& A' & \rightarrow
     & 0 \
  \end{array}
  $$
   be a commutative diagram such that in the horizontal short exact
   sequences epimorphisms have completely positive and contractive sections.
   Then diagram
   $$
   \begin{array}{ccc}
     E_n(I) & \stackrel{\delta}{\rightarrow} & E_{n-1}(A) \\
    \; \;\;\downarrow ^{E_n(f_I)}&  & \;\;\;\downarrow ^{E_n(f_A)}\\
     E_n(I') & \stackrel{\delta '}{\rightarrow}& E_{n-1}(A') \\
   \end{array}
   $$
   commutes.

   ii. The natural two-sided sequence of abelian groups
   $$
   \cdots \rightarrow E_n(A)\rightarrow E_n(B)\rightarrow
   E_n(I)\stackrel{\delta _n}{\rightarrow }E_{n-1}(A)\rightarrow \cdots
   $$
   is exact.
  \item $\mathbb{E}$ is stable. This means that if
  $e_A:A\rightarrow A\otimes \mathcal{K}$ is a homomorphism defined by a map
  $a\mapsto a\otimes p$, here $p$ is a rank one projection in $\mathcal{K}$,
  then $E_n(e_A)$ is an isomorphism.
\end{enumerate}

In this case one has the properties of cohomology theory like to
the above properties of homology theory. We'll list them.

\begin{theorem}
\label{ctnuc}
Let $\mathbb{E}$ be a Cuntz-Bott cohomology theory. Then there are natural
isomorphisms
\begin{equation}
\label{isocun}
E_{n+1}(A)=E_n(A\otimes \mho _{\mathbb{F}})\;\;\;\mathrm{and}\;\;\;\;
E_{n-1}(A)=E_n(A\otimes \Omega _{\mathbb{F}}).
\end{equation}
\end{theorem}

\begin{corollary}
\label{ntiso}
Let $\mathbb{E}$ and $\mathbb{E}'$ be Cuntz-Bott cohomology theories. If
there exists natural isomorphism $\mu _m:E_m\rightarrow E'_m$ for some $m\in Z$,
then there exists natural isomorphism $\mu _n:E_n\rightarrow E'_n$ for all
$n\in Z$.
\end{corollary}

\section{On the Algebraic $K$-theory of $C^{*}$-categories}

Before introducing our view on algebraic $K$-theory
of $C^{*}$-categories, let us make some more comment on the results of
A. Suslin and M. Wodzicki in algebraic $K$-theory. It is well known
fact that any $C^{*}$-algebra has a right (left) bounded approximate unit.
They have, by Proposition 10.2 of \cite{suw}, the factorization property
$(\mathbf{TF})_{\mathrm{right}}$. Thus any $C^{*}$-algebra
possesses property $\mathbf{AH_{Z}}$ and according on Proposition
1.21 of \cite{suw}, one concludes that $C^{*}$-algebras satisfy
excision in algebraic $K$-theory. The results mentioned above leads to
the anew definition of algebraic $K$-theory of  $C^{*}$-algebras,
which is flexible in the comparison of it with the topological $K$-theory
of $C^{*}$-algebras.

Let $A$ be a $C^{*}$-algebra and denote be $A^{+}$ the
$C^{*}$-algebra obtained by adjoining a unit to $A$. If $A$ is
unital denote by $GL_{n}(A)$ the group of invertible elements in
the $C^{*}$-algebra $M_{n}(A)$, and in the non-unital case define
$GL_{n}(A)$ as the group $\ker(GL_{n}(A^{+})\rightarrow
GL_{n}(A))$. Since any $C^{*}$-algebra has a right (left) bounded
approximate unit, it implies well known fact $A=A^{2}$. Thus, by
Corollary 1.13 of \cite{suw} the group of elementary matrices $E(A)$
is a perfect normal subgroup of $GL(A)$ with an abelian guotient
$GL(A)/E(A)$. So one can apply Quillen plus construction to the
classifying space $BGL(A)$. The resulting space denote by $BGL(A)^+$.

The algebraic $K$-theory groups are defined by the following manner:
$$
K_n^a(A)=
  \begin{cases}
     K_0(A)& \text{if}\,\, n=0\\
    \pi _{n-1}(B^{+}(GL(A))  & \text{if} \,\,n\in {\rm N}.
  \end{cases}
$$
and for negative $n$ the group $K_n^a(A)$ (so called Bass
$K$-groups, which sometimes will be denoted by $K_n^B(A)$) is
defined so that the following sequence
$$
K_{1-\,n}^a(A)\rightarrow
K_{1-\,n}^a(A[x,x^{-1}]\rightarrow K_n^a(A)\rightarrow 0
$$
is exact.

Now,according to the results of \cite{suw}, one has the following properties
of algebraic $K$-theory of $C^{*}$-algebras.

\begin{enumerate}
  \item $K_{i}$ is a covariant functor from the category of $C^{*}$-algebras
  and their $*$-homomorphisms into the category of abelian groups for any
  $i\geq 1$;
  \item For every unital $C^{*}$-algebra $R$, which contains a $C^{*}$algebra $A$ as
  a two-sided ideal, the canonical map $K_{*}(A)\rightarrow K_{*}(R,A)$
  an isomorphism;
  \item The natural embedding in the left upper corner
  $A\hookrightarrow M_{n}(A)$ induces, for every natural $n$, an isomorphism
  $K_{*}(A)\simeq K_{*}M_{n}(A)$;
  \item Any extension of $C^{*}$-algebras
  $$
  0\rightarrow I\rightarrow B\rightarrow A\rightarrow 0
  $$
  induces functorial and infinite two-sided long exact sequence of algebraic
  groups
\begin{equation}\label{lsakt}
  \cdots \rightarrow K_{i+1}(A)\rightarrow K_{i}(I)\rightarrow K_{i}(B)
  \rightarrow K_{A}(I)\rightarrow \cdots  \;\;\;\;\;\; (i\in \mathbb{Z}).
\end{equation}
  \end{enumerate}

Using the property (3) and Lemma 2.6.12 in \cite{hig3}, one gets the
following property

\begin{itemize}
  \item ({\it Invariance under inner automorphism}) Let $A$ be a
$C^{*}$-algebra and $u$ be an unitary element in an unital $C^{*}$-algebra
containing $A$ as a closed two-sided ideal. Then the inner automorphism
$ad(u):A\rightarrow A$ induces identity map of $K$-groups.
\end{itemize}

  \subsection{Algebraic $K$-functors of $C^{*}$-categoroids \label{ctnuf}}

 In this subsection we define algebraic $K$-theory of
$C^*$-categoroids in the form which suites our purposes in the sequel.

Let $J$ be a $C^*$-categoroid and $A$ be an additive $C^{*}$-category
containing $J$ as a closed $C^{*}$-ideal.
Let $\mathcal{L}(a)=\hom _A(a,a)$ and $\mathcal{L}(a,J)=\hom _J(a,a)$.
The latter is a closed ideal in the $C^{*}$-algebra $\mathcal{L}(a)$.
A morphism $v:a\rightarrow a'$ in $A$ will be called {\it isometry}
if $v^{*}v=1_a$. Let us write $a\leq a'$ if there is an isometry
$v:a\rightarrow a'$. The relation ''$a\leq a$'' makes the
set of objects into direct system. Any isometry $v:a\rightarrow a'$
defines $*$-homomorphisms of $C^{*}$-algebras
$$
{\rm Ad}(v):\mathcal{L}(a)\rightarrow \mathcal{L}(a')
$$
by the rule $x\mapsto vxv^{*}$. It maps $\mathcal{L}(a,J)$ into
$\mathcal{L}(a',J)$.

Let $v_1:a\rightarrow a'$ and $v_2:a\rightarrow a'$ are
two isometries. Then ${\rm Ad}v_1$ and ${\rm Ad}v_2$ induce
the same homomorphisms
$$
{\rm Ad}_{*}v_1={\rm Ad}_{*}v_2:K_n^a(L(a))\rightarrow
K_n^a(L(a'))
$$
and
$$
{\rm Ad}_{*}v_1={\rm Ad}_{*}v_2:K_n^a(L(a,J))\rightarrow
K_n^a(L(a',J)).
$$
Indeed, let $u$ be an unitary element in an unital $C^{*}$-algebra
containing $A$ as an ideal, then for the inner automorphism
$ad(u):A\rightarrow A$ the homomorphism $K_n^a(ad(u))$ is the identity
map. Therefore, the maps
$$
x\mapsto \left(
\begin{array}{cc}
x & 0 \\
0 & 0
\end{array}
\right) \;\;\;{\rm and\;\;\;}x\mapsto \left(
\begin{array}{cc}
0 & 0 \\
0 & x
\end{array}
\right)
$$
taking $\mathcal{L}(a')$ into $M_2(\mathcal{L}(a'))$,
induces the same isomorphisms after using the functor
$K_n^a$. So, it is enough to show that the maps
$$
x\mapsto \left(
\begin{array}{cc}
v_1xv_1^{*} & 0 \\
0 & 0
\end{array}
\right) \;\;\;{\rm and\;\;\;}x\mapsto \left(
\begin{array}{cc}
0 & 0 \\
0 & v_2xv_2^{*}
\end{array}
\right) ,
$$
which take $\mathcal{L}(a)$ into $M_2(\mathcal{L}(a'))$,
induce by $K_n^a$ the same map. Indeed, the second is obtained
from the first conjugating by the unitary
$$
\left(
\begin{array}{cc}
1-v_1v_1^{*} & v_1v_2^{*} \\
v_2v_1^{*} & 1-v_2v_2^{*}
\end{array}
\right) ,
$$
which is an element of $M_{2}(\mathcal{L}(a'))$.

This discussion shows that the homomorphism
$\nu _{*}^{aa'}=K_n^a(\nu ^{aa'})$
is not depended on choosing an isometry $\nu ^{aa'}:a\rightarrow a'$.
Therefore one has direct system
$\{K_n^a(\mathcal{L}(a)),\nu_{*}^{aa'})\}_{a,a'\in obA}$.

\begin{definition}
\label{wb} Let $J$ be an additive $C^{*}$-categoroid and $A$ be an
additive $C^*$-category containing $J$ as a closed
$C^{*}$-ideal. Define
\begin{equation}
\label{ktg}
\begin{array}{c}
\mathbf{K}_n^a(A,J)=\underrightarrow{\lim }K_n^a(\mathcal{L}(a,J)).
\end{array}
\end{equation}
\end{definition}

\begin{lemma}
\label{wi}
Let $J$ be a $C^{*}$-categoroid considered as a closed ideal in an additive
$C^{*}$-category $A$.
Then $\mathbf{K}_n^a(A,J)=\mathbf{K}_n^a(M(J),J)$.
where $M(J)$ is the multiplier (additive) $C^{*}$-category of $J$.
\end{lemma}

\begin{proof}
Since there exists natural $*$-functor $G:A\rightarrow M(J)$ identity on
$J$, the relation "$a\leq a'$" in $A$ implies "$a\leq a'$" in $M(J)$. This means
that there is natural morphism of direct systems
$$
\{K_n^a(\mathcal{L}(a)),\nu_{*}^{aa'})\}_{a,a'\in obA}\rightarrow
\{K_n^a(\mathcal{L}(a)),\nu_{*}^{aa'})\}_{a,a'\in obM(J)},
$$
which is given by identity homomorphism on each $K_n^a(\mathcal{L}(a))$.
This morphism is cofinal, since if "$a\leq a'$" in $M(J)$ then
"$a\leq a\oplus a'$" in $A$ and "$a'\leq a\oplus a'$" in $M(J)$.
\end{proof}

According to Lemma \ref{wi}, one has the following.

\begin{definition}
\label{waa}Let $J$ be an an additive $C^{*}$-categoroid. Then by definition
$$
\mathbf{K}_n^a(J)=\mathbf{K}_n^a(M(J),J).
$$
\end{definition}

\begin{definition}
\label{lanif}
Let $A$ and $B$ are $C^*$-categories. A $*$-functor $G:A\rightarrow B$
is said to be cofinal if for of a for any object $b\in B$ there exists
an object $a\in A$ and an isometry $s:b\rightarrow G(a)$.
A $C^{*}$-category $A$ is said to be a $cofinal$ sub-category in $B$,
if natural inclusion functor is cofinal.
\end{definition}

The following lemma is trivial but useful in the next part of paper.

\begin{lemma}
\label{gaga}
Let $A'$ be a cofinal full sub-$C^*$-category of an additive $C^*$-category.
Then  $\mathbf{K}_n^a(A')=\mathbf{K}_n^a(A)$.
\end{lemma}
Now, we prove excision property of algebraic $K$-theory on the
category of $C^{*}$-categoroids, which plays major role in this
paper.

\begin{proposition}
\label{utul}Let $A$ and $B$ be additive $C^{*}$-categories and
$J$ be an ideal in $B$ such that the sequence
$$
0\rightarrow J\rightarrow B\rightarrow A\rightarrow 0
$$
is exact. Then two-side sequence of algebraic
$K$-groups
\begin{equation}
\label{exact}
\begin{array}{cc}
...\rightarrow \mathbf{K}_n^a(A)\rightarrow
\mathbf{K}_{n-1}^a(J)\rightarrow \mathbf{K}_{n-1}^a(B)\rightarrow ... &  \\
...\rightarrow \mathbf{K}_0^a(A)\rightarrow ... & ...\rightarrow
\mathbf{K}_{-m}^a(A)\rightarrow \mathbf{K}_{-m-1}^a(J)\rightarrow ...
\end{array}
\end{equation}
is exact $n,m\in {\rm N}$.
\end{proposition}

\begin{proof} Consider exact sequence of $C^{*}$-algebras
$$
0\rightarrow \mathcal{L}(a,J)\rightarrow \mathcal{L}(a)\rightarrow
\mathcal{L}(a)/\mathcal{L}(a,J)\rightarrow 0.
$$
By excision property of algebraic $K$-theory of $C^{*}$-algebras
\cite{suw}, one has two sided long exact sequence
$$
\cdots \rightarrow K_n^a(\mathcal{L}(a)/\mathcal{L}(a,J))\rightarrow
K_{n-1}^a(\mathcal{L}(a,J))\rightarrow K_{n-1}^a(\mathcal{L}(a))\rightarrow
K_{n-1}^a( \mathcal{L}(a)/\mathcal{L}(a,J))\rightarrow \cdots
$$
Using lemma \ref{wi} and definitions \ref{waa} and \ref{wb}, one
can have exactness of (\ref{exact}), because direct limit preserves
exactness.
\end{proof}

\subsection{Comparison of "+" and "Q" Variants of Algebraic $K$-theory
  of $C^*$-categories \label{pmoc}}
  Let $A$\ be a $C^{*}$-algebra with unit. Denote by ${\rm F}(A)$
the additive $C^{*}$-category which has standard Hilbert right
$A$-modules $A^n=A\oplus _{{\rm n-times}}\oplus A$ as the objects
and usual $A$-homomorphisms (which has an adjoint) as the
morphisms. Let ${\rm P}(A)$ be the standard pseudo-abelian
$C^{*}$-category of ${\rm F}(A)$, i.e. ${\rm P}(A)={\rm P(F}(A))$.

If an additive full sub-category $A$ in a pseudo-abelian
category $A'$ is cofinal subcategory then we'll say that $A$ generates $A'$.

Let $a$ be an object in a pseudo-abelian $C^{*}$-category $A$. Put
$a^{\oplus _n}=a\oplus \cdots _{\mathrm{n-times}}\cdots \oplus a,\;\;n\in
\mathbb{N}$. Let $F_{a}$ be full additive subcategory of $A$ which has
$\{a^{\oplus _n}|n=1,...\}$ as the set of objects. Consider a full
sub-$C^{*}$-category $A_a$ consisting all such objects $a'$ in $A$
for which there exists an isometry $s:a'\rightarrow a^{\oplus _n}$, where
$a^{\oplus _n}\in F(a)$. It is clear, that $A_a$ is a pseudo-abelian
$C^*$-category, which is said to be {\it maximal pseudo-abelian sub-}
$C^{*}${\it -category } generated by an object $a\in A$.

\begin{lemma} \label{qmx} Let $A$ be an additive $C^{*}$-category
and $J:A\rightarrow B$ be a $*$-additive functor, where $B$ is a
pseudo-additive $C^*$-category. Then there is a functor
$\mathrm{J}:\mathrm{P}(A)\rightarrow B$ extending $J$. Moreover,
if $J':A\rightarrow B$ is an other functor isomorphic to $J$ then
$\mathrm{J}$ is isomorphic to $\mathrm{J}'$, where latter functor
is an extension of  $J'$. \end{lemma}

\begin{proof} Let $a$ be an object in $A$ and $p_{a}\in
\mathcal{L}(a)$ a projection. Since $B$ is a pseudo-abelian
category, one can choose an object $[p_{a}]$ in $B$ such that
$[1_{a}]=J(a)$ and an isometry $s_{p_{a}}:[p_{a}] \rightarrow
J(a)$ such that $ss^{*}=J(p_{a})$. Define a functor
$\mathrm{J}:\mathrm{P}(A)\rightarrow B$ by the maps $(a,p)\mapsto
[p_{a}]$ and $f\mapsto s_{p'_{a'}}^{*}\cdot J(f)\cdot s_{p_{a}}$,
where $f:(a,p_{a})\rightarrow (a',p'_{a'})$ is a morphism in
${\rm P}(A)$. Simple checking shows that $\mathrm{J}$ satisfies
the requirement of the lemma. To show the second part of the lemma,
let $p\in \mathcal{L}(a)$ be a projection. Then there is an
isometry $\{p\}:(a,p)\rightarrow (a,1)$ induced by projection
$p$. So, one has the isometries $s_{p}:[p]\rightarrow J(a)$ and
$s'_{p}:[p]\rightarrow J'(a)$ such that $\mathrm{J}(a,p)=[p]$,
$\mathrm{J}'(a,p)=[p]'$, and $\mathrm{J}(\{p\})=s_{p}$ and
$\mathrm{I}'(\{p\})=s'_{p}$. It is easily to check that the
collection $\{s'_{p}\cdot \tau _{a}\cdot s_{p}\}$ is the natural
isomorphism of functors from $\mathrm{J}$ to $\mathrm{J}'$, where
$\{\tau _{a}\}$ is a natural isomorphism from $J$ to $J'$.
\end{proof}

\begin{remark} \label{proj} Let $A$ be an unital $C^*$-algebra.
Then Lemma \ref{qmx} implies that ${\rm P}(\mathrm{F}(A))$
(further it is also denoted by ${\rm P}(A)$) is equivalent to the
category ${\cal P}(A)$\ of finitely generated projective right
$A$-modules, where $\mathrm{F}(A)$ is the additive
$C^{*}$-category of standard finitely generated free right
Hilbert $A$-modules. From now on, we consider the pseudo-abelian
$C^{*}$-category ${\rm P}(A)$\ as the substitute of the category
${\cal P}(A)$. \end{remark}

Now, we give an interpretation of Quillens $K$-groups
\cite{quil}, \cite{inas} on which is based our calculations in the
next part of paper.

\begin{definition}
\label{algk}Let $A$ be a pseudo-abelian $C^{*}$-category. Under
$K_n^Q(A)$, $ n\geq 0$ we mean Quillen's $K$-groups relative to
the family of split short exact sequences in $A$; when $A$ is
unital $C^{*}$-algebra, by definition $K_n^Q(A)=K_n^Q({\rm
P}(A))$, $n\geq 0$.
\end{definition}

Let $A$ and $B$ be unital $C^{*}$-algebras and $\varphi
:A\rightarrow B$ be a $*$-homomorphisms (not unital in general).
Then one has a $*$-functor ${\rm P}(\varphi ):{\rm
P}(A)\rightarrow {\rm P}(B)$ defined by the maps $(A^n,p)\mapsto
(B^n,\varphi ^n(p))$ and $(f_{ij})\mapsto (\varphi (f_{ij}))$
where $(f_{ij})$ is $n\times m$-matrix which defines a
$A$-homomorphism from $A^n$ into $A^m$. Therefore we get a functor
${\rm P}$ from the category of unital $C^{*}$-algebras and
$*$-homomorphisms (non-unital in general) into category of
pseudo-abelian $C^{*}$-categories and $*$-functors.

\subsubsection{Construction \label{cons1}} Let $A$ be a pseudo-abelian
$C^{*}$-category and let $s:a'\rightarrow a$ be an isometry in
$A$. There are a $*$-homomorphism $f_s:{\cal L} (a')\rightarrow
{\cal L}(a)$ defined by the map $\alpha \mapsto s\alpha s^{*}$ and
the induced $*$-functor
$$
{\rm P}(f_s):{\rm P(}{\cal L}(a'))\rightarrow {\rm P}({\cal
L}(a)).
$$
If $s_1:a'\rightarrow a$ is an other isometry there exists natural
isomorphism $\upsilon :{\rm P}(f_s)\rightarrow {\rm P}(f_{s_1})$
defined by the following way. For an object of form $(a',p)$ in
${\rm P}({\cal L}(a'))$ let's define $\upsilon
_{(a',p)}:(a,sps^{*})\rightarrow (a,s_1ps_1^{*})$ by equality
$\upsilon _{(a',p)}=s_1ps^{*}$. Since
$$
(sps_1^{*})(s_1ps^{*})=sps^{*}\;\;{\rm and}\;
\;(s_1ps^{*})(sps_1^{*})=s_1ps_1^{*},
$$
$\upsilon _{(a^{\prime },p)}$ is an isomorphism in ${\rm P}({\cal
L}(a))$. In general, for the objects of form $(a^{n'},p)$ we
define $\upsilon _{(a^{n'},p)}=s_1^nps^{n*}$. Since isomorphic
additive functors induce the same homomorphism after using
algebraic $K$-functor, one has
\begin{equation}
\label{isoiso}K_n^a({\rm P}(f_s))=K_n^a({\rm P}(f_{s_1})).
\end{equation}

Let $a$ and $a'$ be objects in $A$. we write $a^{\prime }\leq a$
if there is an isometry $s:a^{\prime }\rightarrow a$. The relation
''$a^{\prime }\leq a$'' makes the set of objects of a
pseudo-abelian category $A$ into a direct system $\{obA,\leq \}$.
Therefore one has correctly defined direct system of abelian
groups.
\begin{equation}
\label{theta1}\Omega _1(A)=\{K_n^a({\cal L}(a)),\kappa _{a^{\prime
}a}^n)\}_{\{obA,\leq \}}
\end{equation}
where $\kappa _{a^{\prime }a}^n$ is the homomorphism $K_n^a({\rm
P}(f_s))$ and by (\ref {isoiso}) it is not depended from the
choosing of an isometry $s:a'\rightarrow a$.

\subsubsection{Construction \label{cons2}} Consider, now, a second direct system of
abelian groups for $A$. Let, as above, $A_a$ be the maximal
pseudo-abelian sub-$C^{*}$-category generated by an object $a\in
A$. It is evident that if there exists isometry $s:a'\rightarrow
a$ then one has natural inclusion additive $*$-functor (not
depended on $s$) $i_{a^{\prime }a}:A_{a^{\prime }}\rightarrow A_a$
and thus we have the direct system $\{A_a,i_{a^{\prime
}a}\}_{(obA,\leq )}$ of the pseudo-abelian $C^{*}$-categories. Let
$\mu _{a'a}^n=K_n(i_{a'a})$. Therefore we have the following
direct system of abelian groups
\begin{equation}
\label{theta2}\Omega _2(A)=\{K_n(A_a),\mu _{a^{\prime
}a}^n\}_{(obA,\leq )},
\end{equation}
which is connected to the direct system $\Omega _1(A)$.

There is a natural isomorphism from the direct system $\Omega
_1(A)$ into the direct system $\Omega _2(A)$. Indeed, consider
natural $*$-functor $\omega _a:{\rm F}({\cal L}(a))\rightarrow A$
which is given by the maps ${\cal L}^n(a)\mapsto a^n$ and
$(f_{ij})\mapsto (f_{ij})$. Then by Lemma \ref{qmx} one can
choose  $\omega _{a}':{\rm P}({\cal L}(a))\rightarrow A$, an
extension of $\omega _{a}$ for every $a\in A$. Elementary
checking shows that $\omega _{a}'$ is equivalence from ${\rm
P}({\cal L}(a))$ onto $A_{a}$.

\begin{proposition} \label{dido}
On the category of pseud-abelian
$C^{*}$-categories and additive $*$-functors the functors
$\mathbf{K}^{a}$ and $\mathbf{K}^{Q}$ are naturally isomorphic.
\end{proposition}

\begin{proof} Let $s:a\rightarrow a'$ be an isometry. Then one
has the isometries $s^{n}:a^{n}\rightarrow a'^{n}$ for any
natural $n$, where $s^{n}=s\oplus\cdots
_{\mathrm{n-times}}\cdots\oplus s$. Define a functor $f
_{s}:\mathrm{P}(a)\rightarrow \mathrm{P}(a')$ by the maps
$(a^n,p)\mapsto (a'^n,s^nps^{*n})$ and $l\mapsto s^mls^{*n}$, where
$l:(a^n,p)\rightarrow (a^m,q)$ is an morphism in $\mathrm{P}(a)$.
We assert that the following diagram
\begin{equation} \label{ckaro}
\begin{array}{ccc}
{\rm P}({\cal L}(a')) & \stackrel{\omega _{a'}}{\rightarrow } & A_{a'} \\
^{f_s}\uparrow \;\; &  & \;\;\cup ^{i_{aa'}}  \\
{\rm P}({\cal L}(a)) & \stackrel{\omega _a}{\rightarrow } & A_a
\end{array}
\end{equation}
is commutative up to isomorphism of functors,
i.e. $\omega _{a'}\cdot f_s\approx i_{a'a}\cdot
\omega _{a}$.
According to Lemma \ref{qmx}, it is enough to construct
an isomorphism
$$
g:\omega _{a'}(f_s(\mathcal{L}(a)))\rightarrow
i_{aa'}(\omega _{a}(\mathcal{L}(a))).
$$
But $\omega _{a'}'(f_s(\mathcal{L}(a)))=\omega_{a'}'(a',ss^{*})$ and
$i_{aa'}(\omega _{a}(\mathcal{L}(a)))=a$. Note that in the
definition of $\omega '$ the object $\omega '((a',ss^{*})=[ss^*]$ is taken
such that there exists isometry $s_{1}:[ss^*]\rightarrow a'$ that
$s_{1}^*s_{1}=1_{ss^*}$ and $s_{1}s_{1}^*=ss^*$. Now, by definition
$g=s^*s_{1}$ which, of course, is an isomorphism from $[ss^*]$ into
$a$. Therefore the latter diagram is commutative up to isomorphism of
functors. Now, apply $K^Q$-functor. Since isomorphic functors have the same
$K$-value, the diagram
\begin{equation}\label{srd}
\begin{array}{ccc}
K_n^Q(\mathcal{L}(a)) & \stackrel{\zeta _n^{a}}{\rightarrow }
 & K_n^Q(A_{a}) \\ ^{v_{a'a}^n}\downarrow
\;\; &  & \;\;\;\downarrow ^{\mu _{a'a}^n} \\
K_n^Q(\mathcal{L}(a')) & \stackrel{\zeta _n^{a'}}{\rightarrow } &
K_n^Q(A_{a'})
\end{array}
\end{equation}
is commutative, where $\zeta _n^a=K^Q(\omega _a)$. Thus we give the natural
isomorphism of direct systems
$\{\zeta _n^a\}:\Omega _{1}\rightarrow \Omega _{2}$.
It is clear, that algebraic ($C^{*}$-algebraic) direct limit
of the direct system $\Omega _2$ is the category $A$. Since
algebraic $K$-functors commute with direct
limits one has
\begin{equation}
\label{mas}
\underrightarrow{\lim}_a K_n^Q(\mathrm{P}(\mathcal{L}(a)))\approx K_{n}^Q(A).
\end{equation}
But $K_n^Q(\mathrm{P}(\mathcal{L}(a)))$ is naturally isomorphic to
$K^a(\mathcal{L}(a))$. Therefore $\mathbf{K}^{a}$ is isomorphic to
$\mathbf{K}^Q$.
\end{proof}

\begin{remark}
\label{nai}Let $A$ be a pseudo-abelian $C^*$ category. Then, for
$n\geq 0$, the groups $\mathbf{K}_n^a(A)$, by Proposition \ref{dido},
is exactly $K_n^Q(A)$. For $n<0$, the groups $\mathbf{K}_n^a(A)$ can
be considered as the generalization of Bass groups, and in this
case this $K$-groups sometimes will be denoted by
$\mathbf{K}_n^B(A)$.
\end{remark}

\section{On the Topological $K$-theory of $C^{*}$-categoroids \label{wash}}

Let $A$ be a $C^{*}$-algebra and $A^{+}$ be the
$C^{*}$-algebra obtained by adjoining a unit to $A$. If $A$ is
unital denote by $GL_{n}(A)$ the topological group of invertible elements
in the $C^{*}$-algebra $M_{n}(A)$, and in the non-unital case define
$GL_{n}(A)$ as the topological subgroup
$$
\ker(GL_{n}(A^{+})\rightarrow GL_{n}(A)).
$$

The topological $K$-theory groups are defined in the following way:
$$
K_n^t(A)=
  \begin{cases}
     K_0(A)& \text{if}\,\, n=0\\
    \pi _{n-1}(GL(A))  & \text{if} \,\,n\in {\rm N}.
  \end{cases}
$$
The following properties of topological $K$-groups.

1. If homomorphisms $f,g:A\rightarrow B$ are homotopic then induced
homomorphisms $K_0^t(f)$ and $K_0^t(g)$ are equal.

2. if
$0\rightarrow I\rightarrow B\rightarrow A\rightarrow 0$
is an exact sequence of $C^*$-algebras then the following sequence
abelian groups
$$
\dots \rightarrow K_n^t(I)\rightarrow K_n^t(B)\rightarrow
K_n^t(A)\stackrel{\delta}{\rightarrow} K_n^t(I)\rightarrow \dots \\
\rightarrow K_0^t(I)\rightarrow K_0^t(B)\rightarrow K_0^t(A)
$$
is exact.

3. Let $A\rightarrow A\otimes \mathcal{K}$
be a homomorphism, defined, for any $C^*$-algebra $A$, by the map
$a\mapsto a\otimes p$, where $p\in \mathcal{K}$ is rank one projection.
Then induced homomorphism is an isomorphism
$K_n^t(A)\simeq K_n^t(A\otimes \mathcal{K})$.

4. Let $\{A_i; f_{ij}\}_I$ be a direct system of $C^*$-algebras and
$\{A;f_i\}$ is the direct limit. Then natural homomorphism
$$
\underrightarrow{\lim} K_n^t(f_i):\underrightarrow{\lim}
K_n^t(A_i)\rightarrow K_n^t(A)
$$
is an isomorphism for any $C^*$-algebra $A$ and natural $n$.

5. By Cuntz-Bott periodicity theorem \cite{cunt}, there are the natural
isomorphisms
$$
K_n^t(A)=
  \begin{cases}
    K_n^t(A\otimes C_0(\mathbb{R})\otimes C_0(\mathbb{R})) &
    \text{in complex case}, \\
      K_n^t(A\otimes C_0(\mathbb{R})\otimes C_0^R(i\mathbb{R}))&
    \text{in real case}.
  \end{cases}
$$
Note that so defined functors has period 2 in the complex $C^*$-algebras
case, and period 8 in the real $C^*$-algebra case.

From the latter property follows that for negative integers $K$-groups may be
defined by the formulas

\begin{itemize}
  \item $K_{-n}^t(A)=K_0^t(A\otimes C_0(\mathbb{R})^{\otimes n})$
for the complex case;
  \item $K_{-n}^t(A)=K_0^t(A\otimes C_0^R(i\mathbb{R}))^{\otimes n}$
  for the real case
\end{itemize}
where $C_0^R(i\mathbb{R}))$ is Cuntz's algebra defined in \cite{cunt}.

From the property 3 and Lemma 2.6.12 in \cite{hig3} immediately follows
that topological (as well as algebraic) $K$-theory has the following
property.
\begin{itemize}
  \item ({\it Invariance under inner automorphism}) Let $A$ be a
$C^{*}$-algebra and $u$ be an unitary element in a unital $C^{*}$-algebra
containing $A$ as a closed two-sided ideal. Then the inner automorphism
$ad(u):A\rightarrow A$ induces identity map of topological $K$-groups.
\end{itemize}

Now, we remark that in the subsection \ref{ctnuf} one can replaces algebraic
$K$-groups by topological $K$-groups then all the results are true. This is
possible since invariance under inner automorphism of algebraic
$K$-theory was used and the same property has topological $K$-theory too.
So we have the following definitions and properties of topological
$K$-theory of $C^*$-categories.

Let $A$ be an additive $C^{*}$-categoroid and $M(A)$ be the
additive $C^{*}$-category centralizers. Let
${\cal {L}}(a)=\hom _{M(A)}(a,a)$ and $A(a)=\hom _A(a,a)$ .
The latter is a closed ideal in the $C^{*}$-algebra ${\cal {L}}(a)$.
Any isometry $v:a\rightarrow a'$ in $M(A)$
defines $*$-homomorphisms of $C^{*}$-algebras
$$
{\rm Ad}(v):A(a)\rightarrow A(a')
$$
by the rule $x\mapsto vxv^{*}$. Thus one has an induced homomorphism
$$
{\rm Ad}_{n}(v):K_n^t(A(a))\rightarrow K_n^t(A(a'))
$$
which isn't depended on choosing of isometry $v:a\rightarrow a'$.
Denote this homomorphism by $\tau _{aa'}$.
We have the direct system of abelian groups
$\{K_n^t(A(a)); \tau _{aa'}\}_{a\in A}$

\begin{definition}
 By definition
\begin{equation}
\label{ktgt}
\begin{array}{c}
{\bf K}_n^t(A)=\underrightarrow{\lim }K_n^t(A(a)).
\end{array}
\end{equation}
\end{definition}

One can easily prove the following.

\begin{lemma}
\label{gagt}Let $A$ be an additive $C^{*}$-category and $A'$ be
a cofinal additive $C^{*}$-subcategory. Then the canonical
additive functor $A'\subset A$ induces an isomorphism
$$
{\bf K}_n^t(A')\approx {\bf K}_n^t(A).
$$
In particular, ${\bf K}_n^t(A)={\bf K}_n^t({\rm P}(A))$.
\end{lemma}

The given bellow is excision property of topological $K$-theory on the
category of $C^{*}$-categoroids.

\begin{proposition}
\label{utut}Let $A$ and $B$ be additive $C^{*}$-categoroids and
$J$ be an ideal in $B$ such that the sequence
$$
0\rightarrow J\rightarrow B\rightarrow A\rightarrow 0
$$
is exact. Then following two-side sequence of topological
$K$-groups
\begin{equation}
\label{exatt}
\begin{array}{cc}
...\rightarrow {\bf K}_n^t(A)\rightarrow {\bf K}_{n-1}^t(J)\rightarrow
{\bf K}_{n-1}^t(B)\rightarrow ... &  \\
...\rightarrow {\bf K}_0^t(A)\rightarrow ... & ...\rightarrow {\bf K}
_{-m}^t(A)\rightarrow {\bf K}_{-m-1}^t(J)\rightarrow ...
\end{array}
\end{equation}
is exact $n,m\in {\rm N}$.
\end{proposition}

Now, we'll give an interpretation of Karoubi's $K$-groups
\cite{kar1}, \cite{kar}  by the functors $K^t$.
The methods, used in the subsection \ref{pmoc}, may be applied.

Since $C^*$-category $A$ is algebraic limit of direct system of
$C^*$-categories $\{A_a, i_{aa'}\}$,
$a\in A$, it is $C^*$-algebraic direct limit too. From the construction of
$C^*$-algebraic direct limit and its property (see subsection \ref{dirl})
implies that the natural homomorphism
$$
\underrightarrow{\lim}(i_{aa'}):\underrightarrow{\lim}K^{-n}(A_a)
\rightarrow K^{-n}(A)
$$
is an isomorphism Karoubi's $K$-groups.

By analogy with Subsection \ref{pmoc}, one can easily proof that $K_n^t$
naturally isomorphic to Karoubi's $K^{-n}$, where $n=0,1,...$

\subsection{Karoubi's Topological $K$-theory of $C^*$-categories}

The purpose of this subsection is to transform some main
results of $K$-theory of Banach categories, introduced by M. Karoubi
in \cite{kar1}, \cite{kar}, to $C^{*}$-categories.

The group $K^0(A)$ of an additive $C^{*}$-category is the
Grothendieck group of the abelian monoid of unitary
isomorphism classes of objects of $A$.
Note that this definition coincides with usual definition because in a
$C^{*}$-category, objects are isomorphic if and only if
they are unitarily isomorphic.
Indeed, if $u:E\rightarrow F$ is isomorphism then $u_0=u\sqrt{(u^{*}u)^{-1}}$
is a unitary isomorphism.

Let $A$ be an additive $C^{*}$-category. The canonical functor induces an
isomorphism $i_{*}:K^0(\tilde A)\rightarrow {\mathbb K}^0(\xi A)$, where the
left-hand $K$-group is the same as in the definition above, and the right one as
in \cite{kar1}, \cite{kar}.

Now we'll give discussion analogous questions for the $K^{-1}$ group
(cf.\ \cite{kar1}, \cite{kar}).
Let $A$ be an additive $C^{*}$-category. Consider the set
of pairs $(E,\alpha)$, where $E\in A$ and $\alpha \in \hom
(E,E)$ is a unitary automorphism.

a). The pairs $(E,\alpha)$ and $(E' ,\alpha ' )$
are said to be \emph{unitarily isomorphic} if there exists a unitary
isomorphism $u:E\rightarrow E'$ such that diagram
$$
\begin{array}{ccc}
E & \stackrel{u}{\rightarrow } & E'  \\
\;\,\downarrow ^\alpha  &  & \;\,\downarrow ^{\alpha '} \\
E & \stackrel{u}{\rightarrow } & E'
\end{array}
$$
is commutative.

b). The pairs $(E,\alpha)$ and $(E,\alpha ')$
are said to be \emph{homotopic} if $\alpha $ and $\alpha ' $ are
homotopic in $\mathrm{Aut}\,E$.

c). A pair $(E,\alpha )$ is said to be \emph{elementary} if it is
homotopic to $(E,1_E)$.

d). The sum is defined by the formula
\[
(E,\alpha)\oplus (E' ,\alpha ')=(E\oplus E',
\alpha \oplus \alpha ').
\]

e). The pairs $(E,\alpha)$ and $(E',\alpha ')$
are said to be \emph{stably isomorphic} if there exist
elementary pairs $(\bar E,\bar e)$
and $(\hat E,\hat e)$, and a unitary isomorphism
\[
(E,\alpha )\oplus (\bar E,\bar e)\simeq (E' ,\alpha ')
\oplus (\hat E,\hat e) .
\]

f). The abelian monoid $K^{-1}(A)$ is defined as the monoid of
classes of stably isomorphic pairs. Denote by $d(E,\alpha)$ the
class of $(E,\alpha)$ in $K^{-1}(A)$.

There are the following relations in $K^{-1}(A)$:

a) $d(E,\alpha )+d(E,\alpha ^{*})=0$\textup{;}

b) If $\alpha $ and $\alpha ' $ are
homotopic unitary isomorphisms,
then $d(E,\alpha )=d(E,\alpha ')$.

c) $d(E,\alpha)+d(E,\beta)=d(E,\beta \alpha)$;

In particular, $K^{-1}(A)$ is an abelian group.

The next proposition is analogous to the
corresponding property of $K^{0}(A)$.

\begin{proposition}
Let $A$ be an additive $C^{*}$-category. The canonical homomorphism
\[
i_{*}:K^{-1}(A)\rightarrow {\mathbb K}^{-1}(A),
\]
defined by $d(E,\alpha )\mapsto d(E,\alpha )$ is an isomorphism.
Here ${\mathbb K}^{-1}(A)$ is Karoubi's group.
\end{proposition}

\begin{proof} $i_{*}$ is an epimorphism: Let $(E,\alpha)$ be a pair
with $\alpha $ an isomorphism. Consider the unitary isomorphism
$\bar \alpha =\alpha \sqrt{\alpha ^{*}\alpha }^{-1}$. It is homotopic
to $\alpha$, because $\alpha ^{*}\alpha $ is homotopic to $1_E$.
We get that $d(E,\alpha )=d(E,\bar \alpha )$. $i$
is a monomorphism: If $i(d(E,\alpha ))=0$, then there exists
elementary $(E' ,e')$ such that $(E\oplus E ,\alpha \oplus e )$
is elementary. Then $(E\oplus E ,\overline{\alpha
\oplus e'})$ is also elementary, that is $(E\oplus E',\alpha
\oplus \bar {e'})$ elementary. This means $d(E,\alpha)=0$.
\end{proof}

Thus the properties of $K^{-1}(A)$ are inherited from the corresponding
properties of ${\mathbb K}^{-1}(A)$. In particular, we get the following:

\begin{theorem}
\label{cvb}Let $A$ be an additive $C^{*}$-category, $\tilde A$ be the
associated pseudoabelian $C^{*}$-category and $i:A\rightarrow \tilde A$ the
canonical additive $*$-functor. Then the induced homomorphism
\begin{equation}
i_{*}:K^{-1}(A)\rightarrow K^{-1}(\tilde A)
\end{equation}
is isomorphism.
\end{theorem}

Let $A$ and $B$ be additive $C^{*}$-categories and $\mathcal{F} :
A\rightarrow B$ be an additive $*$-functor. Denote by
$\Gamma (\mathcal{F})$ the set of triples $(E,F,\alpha)$, where
$E$ and $F$ are objects in $A$, and $\alpha :\mathcal{F}
(E)\rightarrow \mathcal(F)$ is a unitary isomorphism in $B$.

a) Two triples $(E,F,\alpha)$ and $(E',F',\alpha ')$
are \emph{unitarily isomorphic} if there exist unitary
isomorphisms $f:E\rightarrow E'$ and $g:F\rightarrow
F'$ such that the diagram
$$
\begin{array}{ccc}
E & \stackrel{\alpha }{\rightarrow } & F \\
\;\;\;\downarrow ^{f} &  & \;\;\;\downarrow ^{g} \\
E' & \stackrel{\alpha '}{\rightarrow } & F'
\end{array}
$$
is commutative.

b). Two triples $(E,F,\alpha )$ and $(E,F,\alpha ')$
are called \emph{homotopic} if $\alpha$ and $\alpha '$ are
homotopic in the subspace of unitary isomorphisms in
$\hom (E,F)$.

c). The triple $(E,E,1_E)$ is called trivial. A triple
$(E,F,\alpha )$ is said to be
\emph{elementary} if this triple is homotopic to
the trivial triple.

e). The sum of triples is defined by the formula $(E,F,\alpha)\oplus
(E' ,F',\alpha ')=(E\oplus E', F\oplus
F',\alpha \oplus \alpha ')$.

f). Two triples $\sigma =(E,F,\alpha)$ and $\sigma '
=(E',F',\alpha ')$ are \emph{stably unitarily
isomorphic} if there exist elementary pairs $\tau =(\bar E,\bar E,\bar
\alpha)$ and $\tau '=(\bar E ',\bar E'
,\bar \alpha ')$ such that $\sigma \oplus \tau $ and
$\sigma ' \oplus \tau '$ are unitarily isomorphic.

The set $K(\mathcal{F} )$ of stably isomorphic triples is an abelian monoid
with respect to the sum of triples. Denote by $d(E,F,\alpha )$
the class of $(E,F,\alpha )$ in $K(\mathcal{F})$.
The monoid $K(\mathcal{F})$ is an abelian group. Moreover
$
d(E,F,\alpha )+d(F,E,\alpha ^{*})=0.
$
 Note that
$
d(E,F,\alpha )+d(F,E,\alpha ^{*})=d(E\oplus F,F\oplus E,\alpha \oplus \alpha
^{*})
$
The last triple is isomorphic to $(E\oplus F,\beta )$, where
$$
\beta =\left(
\begin{array}{cc}
0 & -\alpha ^{*} \\
\alpha & 0
\end{array}
\right)
$$
which is homotopic to $1_{\mathcal{F}(E)\oplus \mathcal{F}(F)}$ by
$u(t)=\sigma (t)\sqrt{\sigma ^{*}(t)\sigma (t)}$, where
$$
\sigma (t)=\left(
\begin{array}{cc}
1 & -t\alpha ^{*} \\
0 & 1
\end{array}
\right) \left(
\begin{array}{cc}
1 & 0 \\
t\alpha & 1
\end{array}
\right) \left(
\begin{array}{cc}
1 & -t\alpha ^{*} \\
0 & 1
\end{array}
\right)
$$

The following theorem compares our definition of $K(\mathcal{F})$ with the
corresponding one of Karoubi.

\begin{theorem}
The canonical homo\-morphism
$i:K(\mathcal{F})\rightarrow {\mathbb K}(\mathcal{F})$
defined by
$$
d(E,F,\alpha )\mapsto d(E,F,\alpha )
$$
is an isomorphism.
\end{theorem}

\begin{proof} Let $(E,F,\alpha )$ be a triple which defines an
element in ${\mathbb K}(\mathcal{F})$, where $\alpha $ is an isomorphism
(but not unitary isomorphism). Let $\bar \alpha =\alpha
\sqrt{\alpha ^{*}\alpha }$. $\bar \alpha $ is unitary and is
homotopic to $\alpha $ because $\alpha ^{*}\alpha $ is homotopic
to $1_{\mathcal{F}(E)}$. This proves that $i$ is an epimorphism. Now, let
$d(E,F,\alpha)\in K(\mathcal{F})$ defines $0$ in ${\mathbb K}(\mathcal{F})$. This
means, by \cite{kar}, that there exist objects $G$ and $H$ and
isomorphisms (but after polar decomposition we may suppose they
are unitary isomorphisms) $u:E\oplus G\rightarrow H$ and
$v:F\oplus G\rightarrow H$ that $\mathcal{F} (v)(\alpha \oplus 1_{\mathcal{F}
}(G))\mathcal{F}(u^{*})$ is homotopic to $1_{\mathcal{F} (H)}$ (see
\cite{kar}) by a homotopy $h(t)$. Then $\bar
h(t)=h(t)\sqrt{(h^{*}(t)h(t))^{-1}}$ gives homotopy between $(E,F,\alpha
)\oplus (G,G,1_G)$ and $(H,H,1_H)$. This means $d(E,F,\alpha )=0$
in $K(\mathcal{F})$.
\end{proof}

This theorem shows that all properties of $K(\mathcal{F})$ inherited from the
corresponding properties of ${\mathbb K}(\mathcal{F})$. In particular, we'll
get the following results.
(Cf.\ \cite{kar1}, \cite{kar}.)
There are the following relations in $K(\mathcal{F})$:

a) If $\alpha $ and $\alpha ' $ are homotopic, then $d(E,F,\alpha
)=d(E,F,\alpha ')$;

b) $d(E,F,\alpha)+d(F,G,\beta)=d(E,G,\beta\alpha)$.

Let $\mathcal{F} :A\rightarrow B$ be a Serre quasi-surjective additive
$*$-functor. Then

a) if in the definition of $K(\mathcal{F})$ we replace elementary triples by
trivial triples we get the same group.

b) $d(E,F,\alpha )=0$ iff there exist an object $G$ from $A$
and unitary isomorphism $\beta :E\oplus G\rightarrow F\oplus G$ such
that $\mathcal{F}(\beta)=\alpha \oplus 1_{\mathcal{F}(G)}$.

\begin{proposition}
Let $\mathcal{F}:A\rightarrow B$ be a quasi-surjective additive ${*}$-functor.
Then the sequence of abelian groups
\begin{equation}
\label{fgh}K^{-1}(A)\stackrel{f_1}{\longrightarrow}K^{-1}(B)\stackrel{\partial
}{\longrightarrow
}K^0(\mathcal{F})\stackrel{i}{\longrightarrow}K^0(A)\stackrel{
\partial}{\longrightarrow}K^0(B)
\end{equation}
is exact, where $i(d(E,F,\alpha))=d(E)-d(F)$
(for the definition of $\partial$ see \cite{kar}).
In addition, if there  exists a functor $\Psi :B\rightarrow A
$ such that $\mathcal{F} \cdot \Psi \simeq Id_B$, then there exists a split
exact sequence
\begin{equation}
0\rightarrow K^0(\mathcal{F})\stackrel{i}{\rightarrow}K^0(A)\stackrel{j}{
\rightarrow }K^0(B)\rightarrow 0.
\end{equation}
\end{proposition}

Now, we'll discus two examples, which we'll need in the sequel.

1) Recall that an object of $\mathrm{rep}(A,B)$ has the form
$(E,\phi )$, where $E$ is a right Hilbert $B$-module with action of compact group
and $\phi :A\rightarrow \mathcal{L} (E)$ is supposed equivariant. A morphism
from $(E,\phi)$ to $(E', \phi ')$ is by
definition an invariant $B$-homomorphism $f:E\rightarrow E'$
such that $f\phi (a)=\phi ' (a)f$. Note that $\mathrm{rep}(A,B)$ is a
pseudoabelian $C^{*}$-category. To show that $K^i(\mathrm{rep}(A,B))=0$
for all $i\in Z_{2}$,
consider the $\infty $-structure of $\mathrm{rep}(A,B)$ $E^\infty
=E\oplus E\oplus \cdots $, $\alpha ^\infty =\alpha \oplus \alpha
\oplus \cdots $, and $\phi ^\infty (a)=(\phi (a))^\infty $. Let
$$
\infty :\mathrm{rep}(A,B)\rightarrow \mathrm{rep}(A,B)
$$
be the $*$-functor defined by the formula $\infty (E)=E^\infty $,
$\infty (\phi )=\phi ^\infty $, and if $\alpha $ is a
morphism in $\mathrm{rep}(A,B)$, then $\infty (\alpha )=\alpha ^\infty$.
There exists a natural isomorphism
$id_{\mathrm{rep}(A,B)}\oplus \infty \simeq
\infty $. From this it follows that the groups $K^i(\mathrm{rep}(A,B))$
of classes of isomorphic objects of $\mathrm{rep}(A,B)$)
have an automorphism $I$ with property that
$$
id_{K^i(\mathrm{rep}(A,B))}+I=I.
$$
From this fact it follows that $K^i(\mathrm{rep}(A,B))=0$.

2) Consider the canonical quasi-surjective functor
\[
\Theta _{A,B}:\mathrm{rep}(A,B)\rightarrow \mathrm{Cal} (A,B).
\]
Applying the exact sequence (\ref{fgh}) of $K$-groups and result of
example 1, one gets that the canonical homomorphism
\begin{equation}
\label{ghj}\partial :K^{-1}(\mathrm{Cal} (A,B))\rightarrow K^0(\Theta _{A,B})
\end{equation}
is an isomorphism.

\section{Weak Excision \label{auxiso}}

In this section we'll show that the contravariant functors
$$
\mathbf{K}^a_n((Rep(-;B))\;\;\;\mathrm{and}\;\;\;\mathbf{K}^t_n((Rep(-;B))
$$
have weak excision property for all $n\in Z$.

At first one needs following proposition, which easily comes from
the proposition \ref{utul}.

\begin{proposition}
\label{hig} Let $J$ be a $G$-invariant $C^{*}$-ideal in a separable
$G-C^{*}$-algebra $A$, and let
$$
0\rightarrow D(A,J;B)\rightarrow Rep(A,B)\rightarrow
Rep(A,B)/D(A,J;B)\rightarrow 0
$$
be induced exact sequence of $C^{*}$-categoroids. Then for the
algebraic $K$-theory one has the following exact sequences
\begin{equation}
\label{esk}_{...\rightarrow K_n^a(D(A,J;B))\rightarrow
K_n^a(Rep(A,B))\rightarrow \rightarrow
K_n^a(Rep(A,B)/D(A,J;B))\rightarrow K_{n-1}^a(D(A,J;B))\rightarrow
....,}
\end{equation}
respectively for topological $K$-groups
\begin{equation}
\label{ske}_{...\rightarrow K_n^t(D(A,J;B))\rightarrow
K_n^t(Rep(A,B))\rightarrow \rightarrow
K_n^t(Rep(A,B)/D(A,J;B))\rightarrow K_{n+1}^t(D(A,J;B))\rightarrow
....}
\end{equation}
\end{proposition}

Let $0\rightarrow J\rightarrow A\rightarrow A/J\rightarrow 0$ be exact
sequence of $C^*$-algebras such that epimorphism has an equivariant
completely positive and contractive section.
The solution of the problem comes from the showing of
\begin{enumerate}
  \item $\mathbf{K}_n^a(Rep(A;B)/D(A,J;B))\simeq
  \mathbf{K}_n^a(Rep(J;B))$
  \item $\mathbf{K}_n^a(Rep(A/J;B))\simeq \mathbf{K}_n^a(D(A,J;B))$
\end{enumerate}

\subsection{The Isomorphism ${\bf K}_n^a(Rep(A;B)/D(A,J;B))\approx
{\bf K}_n^a(Rep(J;B))$}
Let $(E,\phi )$ be an object in $Rep(A,B)$ and $j:J\rightarrow A$
natural equivariant inclusion.
There is the canonical $*$-functor, induced by the natural inclusion
$j$
\begin{equation}\label{xxx}
{\rm j}:Rep(A;B)\rightarrow Rep(J;B)
\end{equation}
defined by maps $(E,\phi )\mapsto (E,\phi j)$ and $x\mapsto x$.

\begin{proposition}
\label{tipo} The canonical $*$-functor \ref{xxx} maps $D(A,J;B)$
into $D(J,J;B)$ and the induced $*$-functor
\begin{equation}\label{yxy}
\xi:Rep(A;B)/D(A,J;B)\rightarrow Rep(J;B)/D(J,J;B)
\end{equation}
is an isomorphism of $C^{*}$-categories.
\end{proposition}

\begin{proof}(cf. \cite{hig1}) By lemma \ref{print} it is enough to
show that for any object $(E,\phi)$ the $*$-homomorphism
$C^*$-algebras
$$
 \xi _{J,\phi }:D_\phi (A;E;B)/D_\phi (A,J,E;B)\rightarrow
D_{\phi \cdot j}(J,E;B)/D_{\phi \cdot j}(J,J,E;B)
$$
is an $*$-isomorphism. It is easy to show that $\xi _{J,\phi }$ is
a monomorphism. To show that $\xi _{J,\phi }$ is an epimorphism, let
$x\in D_{\phi \cdot j}(J,E;B)$ and let $E_1$ be a $G-C^{*} $-algebra
in $\mathcal{L}(E)$ generated by $\phi (J)\cup \mathcal{K}(E)$; $E_2$ be
the separable $G-C^{*}$-algebra generated by all elements of form
$[x,\phi (y)]$, $y\in J$; and $\mathcal{F}$ be the $G$-invariant separable
linear space generated by $x$ and $\phi (A)$. One has

\begin{itemize}
\item  {$E_1\cdot E_2\subset {\cal K}(E)$,\ \ \ \ because $\phi (b)[\phi
(a),x]\sim [\phi (ba),x]\in {\cal K}(E),\;\;$ $a\in A,\;b\in J$, }

\item  {$[{\cal F},E_1]\subset E_1$, \ \ \ \ \ \ because $[x,\phi
(J)]\subset {\cal K}(E)$ and $[\phi (A),\phi (J)]\subset \phi
(J)$.}
\end{itemize}

From the Kasparov technical theorem follows that there exists
positive $G$-invariant operator $X$ such that
\begin{enumerate}
  \item $X\cdot \phi (J)\subset {\cal K}(E)$;
  \item $(1-X)\cdot [\phi (A),x]\subset {\cal K}(E)$;
  \item $[x,X]\in {\cal K}(E)$.
\end{enumerate}

Since $[(1-X)x,\phi (a)]=(1-X)[x,\phi (a)]-[X,\phi (a)]x$, it
follows from (2) and (3) that $(1-X)x\in D_\phi (A,E;B)$. In
addition, it follows from (2) that $Xx\in D_{\phi \cdot
j}(J,J,E;B)$, and so that the image of $(1-X)x $ in $D_{\phi \cdot
j}(J,E;B)/D_{\phi \cdot j}(J,J,E;B)$ coincides with the image of
$x$.\end{proof}

Now, we prove the following.

\begin{theorem}
\label{topic} Let $A$ be a separable $G-C^{*}$-algebra and $B$ be a
$\sigma $-unital $G-C^{*}$-algebra. Let $J$ be a closed ideal in
$A$. There exists the essential isomorphism
\begin{equation}
\label{qjp}{\bf K}_n^a(Rep(A,B)/D(A,J;B))\approx {\bf
K}_n^a(Rep(J,B))
\end{equation}
\end{theorem}

\begin{proof}
It follows from the Proposition \ref{tipo} that
$$
{\bf K}_n^a(Rep(A;B)/D(A,J;B))\approx {\bf
K}_n^a(Rep(J;B)/D(J,J;B)).
$$
Thus it is enough to show that the natural homomorphism
$$
{\bf K}_{*}^a(Rep(J;B))\rightarrow {\bf
K}_{*}^a(Rep(J;B)/D(J,J;B))
$$
is an isomorphism. From the exact sequence \ref{esk} it follows
that it is enough to show that $K_{*}(D(J,J;B))=0$. By Lemma
\ref{gaga}, we can use the cofinal subcategory $Rep_{H^G_B}(J;B)$,
with the objects of form $(H^G_B,\varphi )$, where $H^G_B$ is
Kasparov's universal Hilbert $B$-module \cite{kas1}. Note that
the canonical isometry
$$
i_1^{H_B}:H^G_B\rightarrow H^G_B\oplus H^G_B
$$
in the first summand is in $D_{\phi ,\phi \oplus
0}(J;H^G_B,H^G_B\oplus H^G_B;B)$ and induces inner homomorphism
$$
{\rm ad(}i_1^{H^G_B}):D_\phi (J,J;H^G_B;B)\rightarrow D_{\phi \oplus
0}(J,J;H^G_B\oplus H^G_B;B).
$$
Consider a sequence of $*$-homomorphisms
\begin{equation}
D_\phi (J,J;H^G_B;B)\rightarrow D_{\phi \oplus \phi }(J,J;H^G_B\oplus
H^G_B;B)\subset D_{\phi \oplus 0}(J,J;H^G_B\oplus H^G_B;B)
\end{equation}
where the inclusion is given by the map $x\mapsto x$ . If the
first arrow is induced by the inclusion
$\iota _1 :H^G_B\rightarrow H^G_B\oplus H^G_B$
in the first summand, then the composition is
$\mathrm{ad}(i_1^{H^G_B})$. If the first arrow is induced by the
inclusion $\iota _2 :H_B\rightarrow H_B\oplus H_B$ in the first
summand, one gives a homomorphism $\lambda ^{H^G_B}$. On the other
hand, the homomorphism $\lambda ^{H^G_B}$ is the composition of the
natural $*$-homomorphisms of $C^*$-algebras
$$
D_\phi (J,J;H^G_B;B)\rightarrow D_0(J,J;H^G_B;B)\rightarrow D_{\phi
\oplus 0}(J,J;H^G_B\oplus H^G_B;B),
$$
given by the maps
$$
x\mapsto x \,\,\mathrm{and}\,\,x\mapsto \left(
\begin{array}{cc}
0 & 0 \\
0 & x
\end{array}
\right).
$$
Remark that $D_0(J,J;H^G_B;B)=M(J\otimes \mathcal{K}_G)$. It is well
known that the latter algebra has trivial algebraic (as well as topological)
$K$-theory groups. If we apply $K$-functors then corresponding
homomorphism of $\lambda ^{H^G_B}$, is zero homomorphisms. Now, let
$ x\in K_{n}^a(D_\phi (J,J;H_B;B))$ represents an element in
$\mathbf{K}_n^a(D(J,J;B))$. Then $x$ and
$\mathbf{K}_n^a(\mathrm{ad}(i_1^{H^G_B})(x))$ represent the same
element. Since
$$
\mathbf{K}_{n}^a(\mathrm{ad}(i_1^{H^G_B}))=\mathbf{K}_{n}^a
(\lambda ^{H^G_B}),
$$
the element represented by $x$ must be zero. Therefore
$\mathbf{K}_{*}(D(J,J;B))=0$.
\end{proof}

\subsection{The Isomorphism
$\Gamma _n:\mathbf{K}_n^a(Rep(A/J;B))\rightarrow \mathbf{K}_n^a(D(A,J;B))$}

Let $A$ and $J$ be as in the last subsection. Let
$p:A\rightarrow A/J$ be the canonical homomorphism which admits an
equivariant completely positive and contractive
section $s:A/J\rightarrow A$. Let $(E,\phi )$ be an object in
$Rep(A;B)$ i.e. there is given an equivariant $*$-homomorphism
$\phi :A\rightarrow \mathcal{L}(E)$. A $*$-homomorphism
$$
\psi =\left(
\begin{array}{cc}
\psi _{11} & \psi _{12} \\
\psi _{21} & \psi _{22}
\end{array}
\right) :A/J\rightarrow \mathcal{L}(E\oplus E'))
$$
will be called $G-s$-$dilation$ for $\phi $ if $\psi _{11}(a)=\phi
(s(a))$, where $E'$ is a right Hilbert $B$-module.

According to generalized Stinespring's theorem there exists a right Hilbert
$B$-module $E'$ and a $G-s$-dilation
$$
\psi =\left(
\begin{array}{cc}
\phi \cdot s & \psi _{12} \\
\psi _{21} & \psi _{22}
\end{array}
\right) :A/J\rightarrow {\cal L}(E\oplus E^{\prime })
$$
for any completely positive and contractive section
$s:A/J\rightarrow A$ \cite{kas1}.

\begin{lemma}
\label{zb}Let $\psi $ be a $s$-dilation for $\phi $. Then\\(1)
$\psi _{12}(a^{*})=\psi _{21}(a)^{*}$;\\(2) for any $a,b\in A$
there exists $j\in J$ such that $\psi _{12}(a)\psi _{21}(b)=\phi
(j)$.\\(3) $\psi _{12}(a)x$ and $x\psi _{21}(a)$ are compact
morphisms for any $a\in A$ and $x\in D_\phi (A,J;B)$.
\end{lemma}

\begin{proof}
{\em The case (1)} is trivial, because $\psi $ is a $
*$-homomorphism. {\em The case (2).} Since $\psi $ is a
$*$-homomorphism, $ \phi \cdot (s(ab)-s(a)\cdot s(b))=\psi
_{12}(a)\cdot \psi _{21}(b)$. But $j=s(ab)-s(a)\cdot s(b)\in J$.
Therefore $\psi _{12}(a)\cdot \psi _{21}(b)=\phi (j)$. {\em The
case} (3). If $x\in D_\phi (A,J;B)$ then $x\phi (j)$ and $\phi
(j)x$ are compact morphisms for any $j\in J$. Then $x\psi
_{12}(a)\cdot \psi _{12}^{*}(a)x^{*}=x\phi (j^{\prime })x^{*}$\
for some $j^{\prime }\in J$. This fact implies that $x\psi
_{12}(a)\cdot \psi _{12}(a^{*})x^{*}$ is compact morphism.
Therefore $x\psi _{12}(a)$ and $\psi _{21}(a)x$ ($=(x^{*}\psi
_{12}(a^{*}))^{*}$) are compact morphisms.
\end{proof}

\begin{lemma}
\label{dila}Let $A$ be separable $C^*$-algebra and $J$ be a
closed ideal such that the projection $p:A\rightarrow A/J$ has a
completely and contractive section $s$. Let $\phi :A/J\rightarrow
\mathcal{L}(E)$ is a $*$-homomorphism and $\psi :A/J\rightarrow
\mathcal{L}(E\oplus E')$ is $s$-dilation for $\phi p$. There exists a
$*$-homomorphism $\varphi :A/J\rightarrow \mathcal{L}(E')$ such that
$$
\psi =\begin{pmatrix}
  \phi & 0 \\
  0 & \varphi
\end{pmatrix}
$$
\end{lemma}

\begin{proof}
A $s$-dilation for $\phi p$ has the form
$$
\begin{pmatrix}
  \phi & \psi _{12}\\
  \psi _{21} & \varphi
\end{pmatrix}.
$$
Since $\psi $ and $\phi$ are $*$-homomorphisms, $\psi _{12}(a)\psi
_{21}(b)=0 $. According to (1) of Lemma \ref{zb}, we have $\psi
_{12}(a)\psi _{12}(a)^*=0 $. Therefore $\psi _{12}(a)=0 $
(similarly, $\psi _{12}(a)=0 $). These facts easily imply that
$\varphi $ is a $*$-homomorphism.
\end{proof}

\begin{lemma}
\label{zs}The map $x=\left(
\begin{array}{cc}
x_{11} & x_{12} \\
x_{21} & x_{22}
\end{array}
\right) \mapsto x^{\prime }=\left(
\begin{array}{ccc}
x_{11} & 0 & x_{12} \\
0 & 0 & 0 \\
x_{21} & 0 & x_{22}
\end{array}
\right) $ defines a $*$-monomorphism
\begin{equation}
\label{ksi}\xi :M_2(D_\phi (A,J,E\oplus E;B))\rightarrow D_{\psi
\cdot p\oplus \phi }(A,J,E\oplus E^{\prime }\oplus E;B).
\end{equation}
\end{lemma}

\begin{proof}
By assumption one has $(\phi (a)\oplus \phi (a))x-x(\phi (a)\oplus
\phi (a))\in {\cal K}(E\oplus E)$, for any $a\in A$, and $(\phi
(b)\oplus \phi (b))x\in {\cal K}(E\oplus E)$, $x(\phi (b)\oplus
\phi (b))\in {\cal K}(E\oplus E)$ for any $b\in J$. It implies
that
$$
\begin{array}{c}
\phi (a)x_{mn}-x_{mn}\phi (a)\in {\cal K}(E),\;\;a\in A,\;\;{\rm
and} \\ \phi (b)x_{mn}\in {\cal K}(E),\;\;x_{mn}\phi (b)\in {\cal
K}(E),\;\;b\in J.
\end{array}
$$
Then $(\psi (p(a))\oplus \phi (a))\cdot x^{\prime }-x^{\prime
}\cdot (\psi
(p(a))\oplus \phi (a))=$%
$$
=\left(
\begin{array}{ccc}
x_{11}\psi _{11}(p(a))-\psi _{11}(p(a))x_{11} & x_{11}\psi
_{12}(p(a)) &
x_{12}\phi (a)-\psi _{11}(p(a))x_{12} \\
\psi _{21}(p(a))x_{11} & 0 & \psi _{21}(p(a))x_{12} \\
x_{21}\psi _{11}(p(a))-\phi (a)x_{21} & x_{21}\psi _{12}(p(a)) &
x_{22}\phi (a)-\phi (a)x_{22}
\end{array}
\right).
$$
As in Lemma \ref{zb}, one has $\psi _{21}(p(a))x_{11}\in {\cal
K}(E,E^{\prime })$, $x_{11}\psi _{12}(p(a))\in {\cal K}(E^{\prime
},E)$, $x_{21}\psi _{12}(p(a))\in {\cal K}(E^{\prime },E)$ and
$\psi _{21}(p(a))x_{12}\in {\cal K}(E,E^{\prime })$. Using the fact
$\phi (a)-\psi _{11}(p(a))\in \phi (J)$, one has
$$
(\psi (p(a))\oplus \phi (a))\cdot x^{\prime }-x^{\prime }\cdot
(\psi (p(a))\oplus \phi (a))\in {\cal K}(E\oplus E^{\prime }\oplus
E),\;\;a\in A.
$$
To show that $(\psi (p(b))\oplus \phi (b))\cdot x^{\prime }$ and
$x^{\prime }\cdot (\psi (p(b))\oplus \phi (b))$ are in ${\cal
K}(E\oplus E^{\prime }\oplus E)$ when $b\in J$, note that $(\psi
(p(b))\oplus \phi (b))\cdot x^{\prime }$ and $x^{\prime }\cdot
(\psi (p(b))\oplus \phi (b))$ are equal to
$$
\left(
\begin{array}{ccc}
0 & 0 & 0 \\
0 & 0 & 0 \\
\phi (b)x_{21} & 0 & \phi (b)x_{22}
\end{array}
\right) \;\;{\rm and}\;\;\left(
\begin{array}{ccc}
0 & 0 & x_{12}\phi (b) \\
0 & 0 & 0 \\
0 & 0 & x_{22}\phi (b)
\end{array}
\right)
$$
respectively. They are compact morphisms because each entries of
matrices are compact.
\end{proof}

Let $A$ be a separable $C^{*}$-algebra, $J$ be a closed ideal in
$A$ and $p:A\rightarrow A/J$ be the canonical $*$-homomorphism.
Let $D^{(p)}(A,J;B)$ be a full $C^{*}$-sub-categoroid in
$D(A,J;B)$ which has all pair of the form $(E,\phi \cdot p)$ as
objects, where a pair $(E,\phi )$ is an object in $D(A/J;B)$.

Consider a $*$-functoroid $\Gamma ^{\prime
}:Rep(A/J;B))\rightarrow D^{(p)}(A,J;B)$ defined by the following
rules:

\begin{enumerate}
\item  if $(E,\phi )$ is an object in $Rep(A/J;B)$, then the corresponding
object in $D^{(p)}(A,J;B)$ is the object $(E,\phi \cdot p)$;

\item  if $x:(E,\phi )\rightarrow (E^{\prime },\phi ^{\prime })$ is a
morphism in $Rep(A/J;B)$, the corresponding morphism is $x:(E,\phi
\cdot p)\rightarrow (E^{\prime },\phi ^{\prime }\cdot p)$.
\end{enumerate}

Let a functoroid $\Gamma :Rep(A/J;B))\rightarrow
D(A,J;B)$ be the composition of $\Gamma ^{\prime }$ with the natural
$*$-inclusion $\varepsilon :D^{(p)}(A,J;B)\subset D(A,J;B)$.

\begin{lemma}
\label{ket}The $*$-functoroid $\Gamma ^{\prime }$ is an
$*$-isomorphism of $C^{*}$-categoroids.
\end{lemma}

\begin{proof} Since $Rep(A/J;B))$ and $D^{(p)}(A,J;B)$ are additive
$C^{*}$-categoroids, it is enough, by Lemma \ref{print}, to show
that for any object $(E,\phi )$ the induced $*$-homomorphism
$$
\Gamma _{(E,\phi )}^{\prime }:D_\phi (A/J;E;B)\rightarrow D_{\phi
\cdot p}(A,J;E;B)
$$
is an $*$-isomorphism. Indeed, if $x\phi (a^{\prime })-\phi
(a^{\prime })x\in {\cal K}(E)$, $\forall a^{\prime }\in A/J$, then
$x\phi (p(a))-\phi (p(a))x\in {\cal K}(E)$, $\forall a\in A$ and
$x\phi (p(j))=\phi (p(j))x=0\in {\cal K}(E)$, $\forall j\in J$.
Conversely is more trivial.
\end{proof}

In the theorem bellow we'll need a $*$-functoroid
$$
\epsilon :D(A,J;B)\rightarrow D^{(p)}(A,J;B)
$$
defined in the following way. For any object $(E,\phi )$ choose an
object $(E\oplus E^{\prime },\psi ^\phi \cdot p)$ such that $\psi
^\phi $ must be a $s$-dilation for $\phi $. This is possible,
since $p$ has completely positive and contractive section $s$, and
by Stinespring's theorem there exists a $s$-dilation for any
completely positive and contractive section. If $x:(E,\phi
)\rightarrow (E^{\prime },\phi ^{\prime })$ is a morphism in
$D(A,J;B)$ then
$$
\left(
\begin{array}{cc}
x & 0 \\
0 & 0
\end{array}
\right) :(E,\psi ^\phi \cdot p)\rightarrow (E^{\prime },\psi
^{\phi ^{\prime }}\cdot p)
$$
is a morphism in $D^{(p)}(A,J;B)$. Indeed,
\begin{equation}
\label{mate}
\begin{array}{c}
\left(
\begin{array}{cc}
x & 0 \\
0 & 0
\end{array}
\right) \left(
\begin{array}{cc}
\psi _{11}^\phi (p(a)) & \psi _{12}^\phi (p(a)) \\
\psi _{21}^\phi (p(a)) & \psi _{22}^\phi (p(a))
\end{array}
\right) -\left(
\begin{array}{cc}
\psi _{11}^{\phi ^{\prime }}(p(a)) & \psi _{12}^{\phi ^{\prime }}(p(a)) \\
\psi _{21}^{\phi ^{\prime }}(p(a)) & \psi _{22}^{\phi ^{\prime
}}(p(a))
\end{array}
\right) \left(
\begin{array}{cc}
x & 0 \\
0 & 0
\end{array}
\right)  \\
=\left(
\begin{array}{cc}
x\psi _{11}^\phi (p(a))-\psi _{11}^{\phi ^{\prime }}(p(a))x &
x\psi
_{12}^\phi (p(a)) \\
-\psi _{21}^{\phi ^{\prime }}(p(a))x & 0
\end{array}
\right) .
\end{array}
\end{equation}
Since
$$
x\psi _{11}^\phi (p(a))-\psi _{11}^{\phi ^{\prime }}(p(a))x=x\phi
(s\cdot p(a))-\phi (s\cdot p(a))x=(x\phi (a)-\phi (a)x)+x\phi
(j)-\phi (j)x,
$$
$\forall a\in A$, where $j=a-s\cdot p(a)\in J$, one has that
$x\psi _{11}^\phi (p(a))-\psi _{11}^{\phi ^{\prime }}(p(a))x$ is a
compact morphism. The morphisms $x\psi _{12}^\phi (p(a))$ and
$\psi _{21}^{\phi ^{\prime }}(p(a))x$ are, by Lemma \ref{zb}, compact
morphism too. Therefore (\ref{mate}) is a compact morphism. Besides,
$$
\left(
\begin{array}{cc}
x & 0 \\
0 & 0
\end{array}
\right) \left(
\begin{array}{cc}
\psi _{11}^\phi (p(j)) & \psi _{12}^\phi (p(j)) \\
\psi _{21}^\phi (p(j)) & \psi _{22}^\phi (p(j))
\end{array}
\right) =\left(
\begin{array}{cc}
\psi _{11}^\phi (p(j)) & \psi _{12}^\phi (p(j)) \\
\psi _{21}^\phi (p(j)) & \psi _{22}^\phi (p(j))
\end{array}
\right) \left(
\begin{array}{cc}
x & 0 \\
0 & 0
\end{array}
\right) =0.
$$
Thus $*$-functoroid $\epsilon $ is correctly defined.

\begin{theorem}
\label{nichbisi} Let $A$ be a separable $C^{*}$-algebra, $J$ be a
closed ideal in $A$ such that canonical $*$-homomorphism
$p:A\rightarrow A/J$ has a completely and contractive section.
Then the functoroid $\Gamma $ induces an isomorphism
$$
\Gamma _n:{\bf K}_n^a(Rep(A/J;B))\rightarrow {\bf K}_n^a(D(A,J;B))
$$
\end{theorem}

\begin{proof}According to Lemma \ref{ket}, it is enough to show that
$\varepsilon :D^{(p)}(A,J;B)\subset D(A,J;B)$ induces the
isomorphism
$$
\varepsilon _n:{\bf K}_n^a(D^{(p)}(A,J;B))\rightarrow {\bf
K}_n^a(D(A,J;B)).
$$
Consider homomorphism
$$
\epsilon _n:{\bf K}_n^a(D(A,J;B))\rightarrow {\bf
K}_n^a(D^{(p)}(A,J;B)).
$$
We assert, that $\epsilon _n\varepsilon _n$ and $\varepsilon
_n\epsilon _n$ are the identity homomorphisms. Let consider the first
case. Let $(E,\phi \cdot p)$ be an object in $D^{(p)}(A,J;B)$. Then
functor $\epsilon \varepsilon $ sends it in the object of form
$(E\oplus E^{\prime },\psi \cdot p)$ where $\psi $ is $s$-dilation
of $\phi \cdot p$. According to Lemma \ref{dila} $\psi =\left(
\begin{array}{cc}
\phi  & 0 \\
0 & \varphi
\end{array}
\right) $ where $\varphi $ is $*$-homomorphism from $A$ into ${\cal
L} (E^{\prime })$. If $x\in D_{\phi \cdot p}(A,J;E;B)$ then
$\epsilon \varepsilon (x)=\left(
\begin{array}{cc}
x & 0 \\
0 & 0
\end{array}
\right) \in D_{\psi p}(A,J;E\oplus E^{\prime };B)$. Let
$i_1:E\rightarrow E\oplus E^{\prime }$ be the inclusion in the
first summand. We assert, that $i_1\phi (p(a))=\psi (p(a))i_1$ and
$\epsilon \varepsilon (x)=i_1xi_1^{*}$ for
$\forall a\in A$. Indeed,
$i_1(\phi (p(a))(\xi ))=\phi (p(a))(\xi )\oplus 0$ and
$$
\psi (p(a))i_1(\xi )=\psi (p(a))(\xi \oplus 0)=\phi (p(a))(\xi )\oplus 0,
\;\;\;\xi \in E.
$$
Now, let $k\in K_n^a(D_{\phi \cdot p}(A,J;E;B))$ represents an element
$\{k\}$ in ${\bf K}_n^a(D^{(p)}(A,J;E;B))$. Then
$\epsilon _n\varepsilon _n(\{k\})=\{\epsilon _n\varepsilon _n(k)\}$.
Since $\epsilon \varepsilon =ad(i_1)$ one has
$\{\epsilon _n\varepsilon _n(k)\}=\{k\}$. Therefore
$\epsilon _n\varepsilon _n$ is identity homomorphism.
Now, we show that homomorphism $\varepsilon _n\epsilon _n$ is the identity
homomorphism. Recall that restriction of $\varepsilon \cdot \epsilon $ on
the object $(E,\phi )$ induces the $*$-homomorphism
$$
D_\phi (A,J;E;B)\rightarrow D_{\psi p}(A,J;E\oplus E^{\prime }B)
$$
given by the map $x\mapsto \left(
\begin{array}{cc}
x & 0 \\
0 & 0
\end{array}
\right) $. Consider a homomorphism $\vartheta$ which is the composition
$$
D_\phi (A,J;E;B)\rightarrow D_{\psi p}(A,J;E\oplus E^{\prime
}B)\rightarrow D_{\psi p\oplus \phi }(A,J;E\oplus E^{\prime
}\oplus E;B),
$$
where the first is induced by $\varepsilon \epsilon $ and the
second arrow is induced by the isometry in the first two summands.
The $*$-homomorphism $\vartheta :D_\phi (A,J;B)\rightarrow D_{\psi
\cdot p\oplus \phi }(A,J;B)$
is defined by
$$
x\mapsto \left(
\begin{array}{ccc}
x & 0 & 0 \\
0 & 0 & 0 \\
0 & 0 & 0
\end{array}
\right) .
$$
Note that this homomorphism has the decomposition:
$$
D_\phi (A,J;B)\stackrel{i_1}{\rightarrow }M_2(D_\phi (A,J;B))
\stackrel{\xi }{\rightarrow }D_{\psi \cdot p\oplus \phi }(A,J;B)
$$
where $i_1$ is given by the map $x\mapsto \left(
\begin{array}{cc}
x & 0 \\
0 & 0
\end{array}
\right) $ and second arrow is (\ref{ksi}). Consider $\eta $ the
composition of the sequence of $*$-homomorphisms
$$
D_\phi (A,J;B)\stackrel{i_2}{\rightarrow }M_2(D_\phi (A,J;B))
\stackrel{\xi }{\rightarrow }D_{\psi \cdot p\oplus \phi }(A,J;B)
$$
where $i_2$ is given by the correspondence $x\mapsto \left(
\begin{array}{cc}
0 & 0 \\
0 & x
\end{array}
\right) $. The homomorphisms $\vartheta$ and $\eta$ induces  the
same homomorphisms of $K$-groups, since $i_1$ and $i_2$ induce the
same homomorphisms of $K$-groups. Let $k\in K_n^a(D_{\phi \cdot
}(A,J;E;B))$ defines an element $\{k\}\in {\bf K}_n^a(D(A,J;B))$.
Let $k'$ be the image of $k$ relative to the homomorphism induced by
the $*$-homomorphism $\vartheta $. Of course, $\{\varepsilon
_n\epsilon _n(k)\}=\{k'\}$. On the other hand, $k'$ coincides with the
image of $k$ relative to the homomorphism induced by the
$*$-homomorphism $\eta $. Since $\eta $ is defined by the
isometry of $E$ in the third summand of $E\oplus E^{\prime
}\oplus E$, one has $\{k\}=\{k^{\prime }\}=\{\varepsilon
_n\epsilon _n(k)\}$. This means that $\varepsilon _n\epsilon _n$
is the identity homomorphism.
\end{proof}

Now, we are ready to prove the following excision property, which plays
one of the main role in this article.

\begin{theorem}
\label{excision}Let $B$ be a $\sigma $-unital $C^{*}$-algebra. If
in an exact sequence of separable $C^{*}$-algebras $0\rightarrow
I\rightarrow A\stackrel{p}{\rightarrow }A/I\rightarrow 0$ $p$
admits a completely positive and contractive section. Then there
exist two-sided long exact sequence
\begin{equation}
\label{exalg} _{...\rightarrow {\bf K}_n^t(Rep(A,B))\rightarrow
{\bf K}_n^a(Rep(J,B))\rightarrow {\bf
K}_{n-1}^a(Rep(A/J,B))\rightarrow {\bf K}
_{n-1}^a(Rep(A,B))\rightarrow }...
\end{equation}
\end{theorem}

\begin{proof}Consider short exact sequence of $C^{*}$-categoroids
$$
_{0\rightarrow D(A,J;B)\rightarrow Rep(A,B)\rightarrow
Rep(A,B)/D(A,J;B)\rightarrow 0.}
$$
We have, by Proposition \ref{utul}, the two-sided long exact
sequence
$$
_{..\rightarrow {\bf K}_n^a(Rep(A,B))\rightarrow {\bf K}
_n^a(Rep(A,B)/D(A,J;B))\stackrel{\partial }{\rightarrow }{\bf K}
_{n-1}^a(D(A,J;B))\rightarrow {\bf K}_{n-1}^a(Rep(A,B))\rightarrow
}
$$

According to Theorem \ref{topic} and Theorem \ref{nichbisi}, one has
\begin{equation}\label{}
\mathbf{K}_n^a(Rep(A;B)/D(A,J;B))\approx \mathbf{K}_n^a(Rep(J;B))
\end{equation}
and
\begin{equation}\label{}
\mathbf{K}_n^a(Rep(A/J;B))\approx \mathbf{K}_n^a(D(A,J;B).
\end{equation}
 It gives us the two-sided long exact sequence (\ref
{exalg}).
\end{proof}

 One has also similar sequence for the case of
topological $K$-theory:
\begin{equation}
\label{exstop}._{..\rightarrow {\bf K}_n^t(Rep(A,B))\rightarrow
{\bf K} _n^t(Rep(J,B))\rightarrow {\bf
K}_{n-1}^t(Rep(A/J,B))\rightarrow {\bf K}
_{n-1}^t(Rep(A,B))\rightarrow ...}
\end{equation}

\section{The Isomorphisms of Algebraic, Topological and Kasparov
$KK$-groups \label{isoatk}}

In this section we'll turn to the our main problem mentioned in the
introduction.

Define algebraic bivariant $KK$-groups by
\begin{equation}\label{kkalg}
KK^a_n(-;B)=
  \begin{cases}
     \mathbf{K}_{n+1}^Q(\mathrm{Rep}(-;B)) & \text{if}, \;\;n\geq -1 \\
     \mathbf{K}_{n+1}n^B(\mathrm{Rep}(-;B)) & \text{if},\;\;n<-1,
  \end{cases}
\end{equation}
and topological $KK$-groups by
\begin{equation}
KK^t_n(-;B)=\mathbf{K}^{-n-1}_{Kar}(\mathrm{Rep}(-;B))
\end{equation}
where $\mathbf{K}_n^Q$, $\mathbf{K}_n^B$ and $\mathbf{K}_{kar}^n$
are Quillen, Bass and Karoubi $K$-groups respectively. We'll write
by $KK_n(A,B)$ Kasparov's group $KK^{-n}(A,B)$. Our main goal is
proof of isomorphism of families mentioned above.

In the first place we proof isomorphism in the fix dimension.

\subsection{On the Isomorphism $\mathbf{K}^a_0(\mathrm{Rep}(-;B))\simeq
KK_{-1}(-;B)$} Consider a triple $(\varphi ,E;p)$, where $E$ is
trivially graded countable generated right $B$-module, $\varphi
:A\rightarrow \mathcal{L}_B(E)$ is a $*$-homomorphism and $p\in
\mathcal{L}_B(E)$, so that
\begin{equation}
\label{kas}
\begin{array}{c}
p\varphi (a)-\varphi (a)p\in \mathcal{K}_B(E), \\ (p^{*}-p)\varphi
(a)\in \mathcal{K}_B(E),\;\;(p^{2}-p)\varphi (a)\in
\mathcal{K}_B(E),\;\;\forall a\in A.
\end{array}
\end{equation}
Such a triple is called {\em Kasparov-Fredholm} $A,B$-{\em
module}. If all left parts in \ref{kas} are zero, then such a
triple is said to be degenerate.

Define sum Kasparov-Fredholm $A,B$-modules by the formula
$$
(\varphi ,E;p)\oplus (\varphi ',E';p')=(\varphi \oplus \varphi
',E\oplus E';p\oplus p').
$$
Consider the equivalence relations:
\begin{itemize}
  \item ({\em Unitary isomorphism}) $A,B$-modules $(\varphi ,E;p)$ and
$(\varphi ',E';p')$ will be said to be unitary isomorphic if there
exists unitary isomorphism $u:E\rightarrow E'$ such that
$$
u\varphi (a)u^{*}=\varphi '(a),\;\;upu^{*}=p',\;\;\forall a\in A.
$$
  \item ({\em Homology}) $A,B$-modules $(\varphi ,E;p)$ and
$(\varphi ',E;p')$ will be said to be homological if
$$
p'\varphi '(a)-p\varphi (a)\in \mathcal{K}_B(E),\;\;\forall a\in
A.
$$
\end{itemize} Simple checking shows that the equivalence relations,
defined above, are well behaved relative to the sum.

Let $\mathcal{E}^1(A,B)$ be a abelian monoid of classes of
equivalence $A,B$-modules generated by the unitary isomorphism and
homology. Denote by $\mathcal{D}^1(A,B)$ a sub-monoid of
$\mathcal{E}^1(A,B)$ consisting of only those classes which are
classes of all degenerate triples. By definition
$$
E^1(A,B)=\mathcal{E}^1(A,B)/\mathcal{D}^1(A,B)
$$

Using the Kasparov stabilization theorem, one can to show that
definition of $E^1(A,B)$ coincides with Kasparov's original
definition of $E^1(A,B)$. Therefore one sees that the last monoid
(group) may be considered as a model of $KK^1(A,B)$ by lemma 2 of
section 7 of \cite{kas2}.

Recall, that the objects in $\mathrm{Rep}(A,B)$, by definition,
has the form $(\varphi ,E;p)$, where $p:(\varphi ,E)\rightarrow
(\varphi ,E)$ is projection in the category $Rep(A;B)$. More
precisely,
$$
\varphi (a)p-p\varphi (a)\in
\mathcal{K}_B(E),\;\;p^{*}=p,\;\;p^2=p.
$$
An unitary isomorphism $s:(\varphi ,E;p)\rightarrow (\psi ,E',q)$
is partial isomorphism $s:E\rightarrow E'$ such that
$$
s\varphi (a)-\psi (a)s\in
\mathcal{K}_B(E,E'),\;\;s^{*}s=p,\;\;ss^{*}=q.
$$
Let $\mathcal{E}_\pi ^1(A,B)$ be abelian monoid of unitary
isomorphic objects in $\mathrm{Rep(A;B)}$. Thus Grothendieck group
of $\mathcal{E}_\pi ^1(A,B)$ is exactly $K_0(\mathrm{Rep}(A,B))$.

Kasparov-Fredholm $A,B$-modules $(\varphi ,E;p)$ and $(\varphi
',E';p')$ will be said to be strong unitary isomorphic if there
exists unitary isomorphism $u:E\rightarrow E'$ such that
$$
u\varphi (a)u^{*}-\varphi '(a)\in
\mathcal{K}_B(E'),\;\;\;\;p'=upu^{*},\;\;\forall a\in A.
$$
Let $s:(\varphi ,E;p)\rightarrow (\varphi ',E';p')$ be an unitary
isomorphism  in $\mathrm{Rep(A;B)}$. Then
$$
\bar s:\left(
\begin{array}{ccc}
\varphi \oplus \psi & E\oplus E & \bar p
\end{array}
\right) \rightarrow \left(
\begin{array}{ccc}
\psi \oplus \varphi & E\oplus E & \bar p'
\end{array}
\right)
$$
is strong isomorphism, where
$$
\bar s=\left(
\begin{array}{cc}
s & 1-ss^{*} \\
1-s^{*}s & s
\end{array}
\right),\;\;\; \bar p=\left(
\begin{array}{cc}
p & 0 \\
0 & 1
\end{array}
\right)\;\;\;\mathrm{and}\;\;\; \bar p'=\left(
\begin{array}{cc}
1 & 0 \\
0 & p'
\end{array}
\right).
$$
Other hand, simple checking shows that strong unitary isomorphic
Kasparov-Fredholm $A,B$-modules is contained in the equivalence
generated by unitary isomorphism and homology. This means that one
has correctly defined homomorphism
$$
\lambda _1:K^0(\mathrm{Rep}(A,B))\rightarrow E^1(A,B),
$$
given by the map $[(\varphi ,E;p)]\mapsto \{(\varphi ,E;p)\}$.

Let $\mathrm{Rep}(A,B)/\mathrm{D}(A,A;B))$ be pseudo-abelian
category of the category $Rep(A,B)/D(A,A;B))$. Let $(\varphi
;E;p)$ be a Kasparov-Fredholm $A,B$-module. The $p$ defines
projector $\dot{p}$ in the category $Rep(A,B)/D(A,A;B))$. Thus the
triple $(\varphi ;E;\dot{p})$ is an object in
$\mathrm{Rep}(A,B)/\mathrm{D}(A,A;B))$. We want definition of a
homomorphism
$$
\mu :E^1(A,B)\rightarrow K_0(\mathrm{Rep}(A,B)/\mathrm{D}(A,A;B))
$$
by the map $(\varphi ;E;p)\mapsto (\varphi ;E;\dot{p})$. We'll
show that this is correct.

We recall definition of operatorial homotopy:
\begin{itemize}
  \item ({\em Operatorial homotopy}) $A,B$-module $(\varphi ,E;p)$ is
operatorial homotopic to a triple $(\varphi ,E;p')$ if there exist
a continuous map $p_t:[0;1]\rightarrow \mathcal{L}_B(E)$ such that
$(\varphi ,E;p_t)$ is $A,B$-module for any $t\in [0;1]$.
\end{itemize}

If $(\varphi ,E;p)$ is homological to $(\psi ,E;q)$, then
$(\varphi ,E;p)\oplus (\psi ,E;0)$ is operatorial homotopic to
$(\varphi ,E;0)\oplus (\psi ,E;q)$. Indeed, Desired homotopy is
defined by the formula
$$
\left( \left(
\begin{array}{cc}
\varphi & 0 \\
0 & \psi
\end{array}
\right) ,E\oplus E,\frac 1{1+t^2}\left(
\begin{array}{cc}
p & tpq \\
tqp & t^2q
\end{array}
\right) \right) ,\;\;\;t\in [0;\infty ].
$$
(cf. section 7 in \cite{kas2}). Thus the projections $\dot{p\oplus
0}$ and $\dot{0\oplus q}$ are homotopic. Then, using Lemma 4
section 6 in \cite{kas2}, one concludes that the objects $(\varphi
,E;\dot{p})\oplus (\psi ,E;\dot{0})$ and $(\varphi
,E;\dot{0})\oplus (\psi ,E;\dot{q})$ are unitary isomorphic
objects in $\mathrm{Rep}(A,B)/\mathrm{D}(A,A;B)$. Let $(\varphi
,E;p)$ is unitary isomorphic to $(\psi ,E;q)$. Then $(\varphi
,E;\dot{p})$ is isomorphic to $(\psi ,E;\dot{q})$ in the category
$\mathrm{Rep}(A,B)/\mathrm{D}(A,A;B)$. Therefore $\mu $ is
correctly defined.

We are just ready to prove the following theorem.

\begin{theorem}
\label{iginol}The natural homomorphism
$$
\lambda _1:K^0(\mathrm{Rep}(A,B))\rightarrow E^1(A,B)
$$
is an isomorphism.
\end{theorem}

\begin{proof} Of course, $ \lambda _1$ is epimorphism. Indeed,
let $(\varphi ,E;p)$ be Kasparov-Fredholm $A,B$-module. Applying
analogue to Lemmas 17.4.2-17.4.3 in \cite{bla}, one can suppose
that $p^*=p$ and $||p||\leq 1$. Then it is equivalent to $(\varphi
\oplus 0,E\oplus E;p')$, where
$$
p'=\begin{pmatrix}
  p & \sqrt{p-p^2}\\
  \sqrt{p-p^2} & 1-p
\end{pmatrix}.
$$
Simple checking shows that $p'$ is a projection and $(\varphi
\oplus 0,E\oplus E;p')$ is an object in $\mathrm{Rep}(A,B)$. To
show that $\lambda _1$ is monomorphism, consider commutative
diagram
$$
\begin{array}{ccc}
K^0(\mathrm{Rep}(A,B)) & \stackrel{\lambda _1}{\rightarrow } & E^1(A,B) \\
        \parallel   &  & \downarrow ^\mu \\
K^0(\mathrm{Rep}(A,B)) & \stackrel{\xi }{\rightarrow } &
K^0(\mathrm{Rep}(A,B)/\mathrm{D}(A,A;B)).
\end{array}
$$
By Theorem \ref{topic}, $\xi $ is an isomorphism. Therefore
$\lambda _1$ is monomorphism.
\end{proof}

\subsection{The Main Theorem}

Now, we present our main result in the following theorem.
\begin{theorem}
\label{corollary} Let $B$ be a $\sigma $-unital $C^{*}$-algebra.
Then the families of functors
\begin{equation}
\label{fmls} \{KK^a_n(-;B)\}_{n\in Z},\;\;\;\{KK^t_n(-;B)\}_{n\in
Z},\;\;\; \{KK_n(-;B)\}_{n\in Z}.
\end{equation}
are naturally isomorphic Cuntz-Bott cohomology theories on the
category of separable $C^*$-algebras and $*$-homomorphisms.
\end{theorem}

\begin{proof}
By Proposition \ref{dido} the functor $KK^a_n(-;B)$ is naturally
isomorphic to $\mathbf{K}^a_{n+1}(\mathrm{Rep}(-;B))$
(respectively for $KK^t_n(-;B)$. See \ref{wash}). But by Theorem
\ref{excision} family
$$
\{ \mathbf{K}^a_n(\mathrm{Rep}(-;B))\}_{n\in Z}.
$$
has weak excision property. The same property has also
$$
\{ \mathbf{K}^t_n(\mathrm{Rep}(-;B))\}_{n\in Z}.
$$
Further, thanks to the result of \cite{cusk}, the family
$$
\{KK_n(-;B)\}_{n\in Z}
$$
has weak excision property. Now, we show:
\begin{itemize}
  \item functors $\mathbf{K}^a_n(\mathrm{Rep}(-;B))$,
  $\mathbf{K}^t_n(\mathrm{Rep}(-;B))$ and $KK_n(-;B)$ have stable
  property, for all $n\in Z$.
\end{itemize}
This fact with the Theorem \ref{excision} implies that all three
families are Cuntz-Bott cohomology theories. Let us go on to show
it.

 Let $p\in {\cal K}$ be a rank one projection and $A$ be a separable
$C^*$-algebra.
The $*$-homomorphism $e_{A}:A\rightarrow A\otimes \mathcal{K}$
defined by the map $a\mapsto a\otimes p$, $\forall a\in A$,
induces a $*$-functor
$$
e_{A}^*:Rep(A\otimes \mathcal{K};B)\rightarrow Rep(A;B).
$$
Now, we construct a $*$-functor
$$
\varepsilon :Rep(A;B)\rightarrow Rep(A\otimes \mathcal{K};B)
$$
which is somehow a right inverse to $e_{A}^*$. Let
$\phi :A\rightarrow \mathcal{L}(E)$
be an object in $Rep(A;B)$. One has the induced $*$-homomorphism
$$
\phi \otimes id_{\mathcal{K}}:A\otimes \mathcal{K}\rightarrow
\mathcal{L}(E\otimes _C\mathcal{H})
$$
defined as the composition
$$
A\otimes \mathcal{K}\rightarrow \mathcal{L}(E)\otimes
\mathcal{L}(\mathcal{H})\rightarrow \mathcal{L}(E\otimes _C\mathcal{H}),
$$
of natural maps. Of course, this is an object in $Rep(A\otimes
\mathcal{K};B)$. Let $f:(E,\phi )\rightarrow (E',\phi ')$ be a
morphism in $Rep(A;B) $, i.e. $f\phi (a)-\phi (a)f\in {\cal
K}(E)$, $\forall a\in A$. Then
$$
\varepsilon (f)=f\otimes p:(E\otimes _C \mathcal{H},\phi \otimes
id_{\mathcal{K}})\rightarrow (E'\otimes _C \mathcal{H},\phi
'\otimes id_{\mathcal{K}})
$$
is a morphism in $Rep(A\otimes \mathcal{K};B)$.

Now we construct an useful isometry $\sigma _E:E\rightarrow
E\otimes _{k}\mathcal{H}$, for any countable generated $B$-module
$E$. Choose $y\in \mathcal{H}$ such that $p(y)=y$ and $||y||=1$
and consider a $B$-homomorphism $\sigma _E$ given by the formula
$x\mapsto x\otimes y$. For a $z\in \mathcal{H}$, there exists
$\lambda _z\in k$ determined uniquely by the equation
$p(z)=\lambda _zy$ (note, that $p$ is rank one projection). Define
$\sigma _E^{*}$ by the correspondence $x\otimes z\mapsto \lambda
_zx$. The $B$-homomorphism $\sigma _E^{*}$ is the adjoint to
$\sigma _E$. Indeed, since $p(y)=y$ and $||y||=1$,
$$
\langle \sigma _E(x);x'\otimes z\rangle =\langle x\otimes
y;x'\otimes z\rangle =\langle x;x'\rangle \cdot \langle
y,pz\rangle =\langle x;\lambda _zx'\rangle =\langle x;\sigma
_E^{*}(x'\otimes z)\rangle \;\; \forall x,x',z\in E.
$$
Since $\sigma _E^{*}\sigma _E(x)= \sigma _E^{*}(x\otimes y)=x$,
one concludes $\sigma _E$ is an isometry.

Since
$$
\sigma _E\phi (a)(z)-((\phi \otimes id_\mathcal{K})e_{A}(a))
\sigma _E(z)=\phi (a)(z)\otimes y-(\phi (a)\otimes p)(z\otimes
y)=0,\;\; \forall z\in E,
$$
the isometry $\sigma _E$ is a morphism from $(E,\phi )$ into
$(E\otimes _C\mathcal{H}, (\phi \otimes id_\mathcal{K})e_{A})$.

Consider restriction of $e^*_{A}\varepsilon $ on the
$D_{\phi}(A;E;B)$. Thus we have the $*$-homomorphism
\begin{equation}
(e^*_{A}\varepsilon )_E:D_{\phi}(A;E;B)\rightarrow D_{(\phi
\otimes id)e_{A}}(A;E\otimes _k \mathcal{H};B)
\end{equation}
which maps $x$ to $x\otimes p$. But
$$
(\sigma _Ex)(z)= (\sigma _E)(x(z))=x(z)\otimes y=((x\otimes
p)\sigma _E)(z),\;\; \forall x\in D_{\phi}(A;E;B),\;\;\forall z\in E.
$$
This means
\begin{equation}
\label{innerel} (e^*_A\varepsilon )_E(x)=\sigma _Ex\sigma _E^*.
\end{equation}
Now, we show the functor $e^*_{A}\varepsilon $ induces the
identity homomorphism of group $\mathbf{K}_n^a(Rep(A,B))$ onto
itself. Indeed, choose an element $r\in \mathbf{K}_n^a(Rep(A,B))$.
By definition of $\mathbf{K}^a_n$-groups the element $r$ is
represented by an element $r_\phi \in K_n^a(D_{\phi}(A;E;B))$.
Then the element $\mathbf{K}^a_n(e^*_{A}\varepsilon )(r)$ is
represented by the element
\begin{equation}
\label{rett}
K^a_n((e^*_{A}\varepsilon)_E)(r_\phi).
\end{equation}
The equation \ref{innerel} implies that the element \ref{rett}
represents the element $r$. This means that $e^*_{A}\varepsilon $
induces identity homomorphism of $\mathbf{K}_n^a(Rep(A,B))$ onto
itself. Thus $\varepsilon _{n}=K_n^a(\varepsilon )$ is a right
inverse of $K_n^a(e^*_{A})$. This means that $K_n^a(e^*_{A})$ is
epimorphism. Thus it is enough to show that $K_n^a(e^*_{A})$ is a
monomorphism. Consider commutative diagram
$$
\begin{array}{ccc}
\mathbf{K}_n^a(Rep(A,B)) &
\stackrel{K_n^a(e^*_{A})}{\longleftarrow}
& \mathbf{K}_n^a(Rep(A^I\otimes \mathcal{K},B))\\
\downarrow ^{\mathbf{K}_n^a(e^A_0)}&  & \downarrow
^{\mathbf{K}_n^a(e^A_0)} \\
\mathbf{K}_n^a(Rep(A^I,B))& \stackrel{\varepsilon
_n}{\longrightarrow} & \mathbf{K}_n^a(Rep(A^I\otimes
\mathcal{K},B)),
\end{array}
$$
where $e^A_t:A^I\rightarrow A$ is the evolution at $t\in I=[0;1]$.
Since $\mathbf{K}_n^a(Rep(-\otimes \mathcal{K},B))$ is split and
stable functor, it is homotopy invariant functor. This fact
implies that right vertical arrow is an isomorphism. Therefore
$K_n^a(e^*_{A})$ is an monomorphism. Therefore the family
$\{\mathbf{K}^a_n(\mathrm{Rep}(-;B))\}_{n\in Z}$ is Cuntz-Bott
cohomology theory.
According to Corollary \ref{ntiso} and Theorem \ref{iginol}, one
concludes that all three theories are isomorphic.
\end{proof}


\end{document}